\documentclass[a4paper]{article}
\usepackage{amssymb}
\usepackage{amsmath}
\usepackage{amsthm}
\usepackage{authblk}
\usepackage[USenglish]{babel}
\usepackage{esint}
\usepackage{graphicx}
\usepackage[colorlinks=true]{hyperref}
\usepackage[latin1]{inputenc}
\newtheorem{theorem}{Theorem}[section]
\newtheorem{lemma}[theorem]{Lemma}
\newtheorem{proposition}[theorem]{Proposition}
\newtheorem{corollary}[theorem]{Corollary}
\newtheorem{remark}[theorem]{Remark}
\newtheorem{definition}[theorem]{Definition}
\newtheorem{example}[theorem]{Example}
\newcommand{\N}{\mathbb N}
\newcommand{\Z}{\mathbb Z}

\newcommand{\R}{\mathbb R}
\newcommand{\C}{\mathbb C}
\renewcommand{\S}{\mathbb S}
\renewcommand{\P}{\mathbb P}
\renewcommand{\H}{\mathbb H}

\newcommand{\T}{\mathbb T}

\newcommand{\mb}{\mathbb}
\newcommand{\mf}{\mathbf}
\newcommand{\mc}{\mathcal}

\newcommand{\mi}{\mathit}
\newcommand{\mr}{\mathrm}

\renewcommand{\a}{\alpha}
\renewcommand{\b}{\beta}
\renewcommand{\d}{\delta}
\newcommand{\D}{\Delta}
\newcommand{\e}{\varepsilon}
\newcommand{\g}{\gamma}
\newcommand{\G}{\Gamma}
\renewcommand{\l}{\lambda}
\renewcommand{\L}{\Lambda}

\renewcommand{\O}{\Omega}
\newcommand{\ph}{\varphi}
\newcommand{\s}{\sigma}
\newcommand{\Si}{\Sigma}
\renewcommand{\t}{\theta}
\newcommand{\z}{\zeta}

\newcommand{\wt}{\widetilde}
\newcommand{\wh}{\widehat}
\newcommand{\ov}{\overline}

\newcommand{\ub}{\underbrace}
\newcommand{\ity}{\infty}
\newcommand{\fin}{\fint}
\newcommand{\fr}{\frac}

\newcommand{\n}{\nabla}

\newcommand{\fa}{\forall}

\newcommand{\es}{\emptyset}
\newcommand{\st}{\star}
\newcommand{\wk}{\rightharpoonup}
\newcommand{\inc}{\hookrightarrow}

\newcommand{\us}{\underset}
\newcommand{\sr}{\stackrel}

\newcommand{\To}{\Rightarrow}

\newcommand{\bs}{\backslash}
\renewcommand{\Cup}{\bigcup}

\newcommand{\sub}{\subset}

\newcommand{\nin}{\not\in}

\newcommand{\ox}{\otimes}

\newcommand{\pl}{\oplus}
\newcommand{\Pl}{\bigoplus}
\newcommand{\seq}{\simeq}

\newcommand{\x}{\times}
\renewcommand{\c}{\circ}
\newcommand{\cd}{\cdot}
\newcommand{\cds}{\cdots}
\newcommand{\dds}{\ddots}
\newcommand{\vds}{\vdots}
\newcommand{\ds}{\dots}
\newcommand{\dis}{\displaystyle}

\newcommand{\lab}{\label}
\newcommand{\eqr}{\eqref}
\newcommand{\ci}{\cite}
\newcommand{\tx}{\text}
\newcommand{\q}{\quad}
\newcommand{\lr}{\left(}
\newcommand{\rr}{\right)}
\newcommand{\ls}{\left[}
\newcommand{\rs}{\right]}
\newcommand{\lb}{\left\{}
\newcommand{\rb}{\right\}}
\newcommand{\lm}{\left|}
\renewcommand{\rm}{\right|}

\newcommand{\ly}{\left.}
\newcommand{\ry}{\right.}

\newcommand{\ti}{\textit}

\newcommand{\cen}{\centering}
\newcommand{\igr}{\includegraphics}
\newcommand{\cpt}{\caption}

\newcommand{\no}{\noindent}

\newcommand{\bst}{\bibliographystyle}
\newcommand{\bdo}{
\begin{document}}
\newcommand{\edo}{\end{document}}
\newcommand{\au}{\author}
\newcommand{\af}{\affil}
\newcommand{\tit}{\title}
\renewcommand{\th}{\thanks}
\newcommand{\ins}{\institute}
\newcommand{\da}{\date}
\newcommand{\mt}{\maketitle}
\newcommand{\bfr}{\begin{frame}}
\newcommand{\efr}{\end{frame}}
\newcommand{\bbl}{\begin{block}}
\newcommand{\ebl}{\end{block}}
\newcommand{\ft}{\frametitle}
\newcommand{\fs}{\framesubtitle}
\newcommand{\pau}{\pause}
\newcommand{\col}{\color}
\newcommand{\ons}{\onslide}
\renewcommand{\sec}{\section}
\newcommand{\sse}{\subsection}
\newcommand{\bab}{\begin{abstract}}
\newcommand{\eab}{\end{abstract}}
\newcommand{\bce}{\begin{center}}
\newcommand{\ece}{\end{center}}
\newcommand{\bHu}{\begin{Huge}}
\newcommand{\eHu}{\end{Huge}}
\newcommand{\bth}{\begin{theorem}}
\renewcommand{\eth}{\end{theorem}}
\newcommand{\ble}{\begin{lemma}}
\newcommand{\ele}{\end{lemma}}
\newcommand{\bpr}{\begin{proposition}}
\newcommand{\epr}{\end{proposition}}
\newcommand{\bco}{\begin{corollary}}
\newcommand{\eco}{\end{corollary}}
\newcommand{\bre}{\begin{remark}}
\newcommand{\ere}{\end{remark}}
\newcommand{\bex}{\begin{example}}
\newcommand{\eex}{\end{example}}
\newcommand{\bde}{\begin{definition}}
\newcommand{\ede}{\end{definition}}
\newcommand{\bpf}{\begin{proof}}
\newcommand{\epf}{\end{proof}}
\newcommand{\bl}{\begin{array}{l}}
\newcommand{\bll}{\begin{array}{ll}}
\newcommand{\ba}{\begin{array}}
\newcommand{\ea}{\end{array}}
\newcommand{\beq}{\begin{equation}}
\newcommand{\eeq}{\end{equation}}
\newcommand{\bea}{\begin{eqnarray}}
\newcommand{\eea}{\end{eqnarray}}
\newcommand{\nn}{\nonumber}
\newcommand{\bey}{\begin{eqnarray*}}
\newcommand{\eey}{\end{eqnarray*}}
\newcommand{\ben}{\begin{enumerate}}
\newcommand{\een}{\end{enumerate}}
\newcommand{\bit}{\begin{itemize}}
\newcommand{\eit}{\end{itemize}}
\renewcommand{\it}{\item}
\newcommand{\bfi}{\begin{figure}}
\newcommand{\efi}{\end{figure}}
\newcommand{\bmi}{\begin{minipage}}
\newcommand{\emi}{\end{minipage}}
\newcommand{\bin}[2]{\left(\genfrac{}{}{0pt}{}{#1}{#2}\right)}
\newcommand{\qm}[1]{``#1''}

\bdo

\tit{Existence and multiplicity result for the singular Toda system}
\au{Luca Battaglia}
\af{S.I.S.S.A., Via Bonomea 265, 34136 Trieste (Italy) - \ti{lbatta@sissa.it}}
\da{}

\mt\

\bab
\no We consider the Toda system on a compact surface $\dis{(\Si,g)}$
$$\lb\bl-\D u_1=2\rho_1\lr\fr{h_1e^{u_1}}{\int_\Si h_1e^{u_1}dV_g}-1\rr-\rho_2\lr\fr{h_2e^{u_2}}{\int_\Si h_2e^{u_2}dV_g}-1\rr-4\pi\sum_{j=1}^J\a_{1j}\lr\d_{p_j}-1\rr\\-\D u_2=2\rho_2\lr\fr{h_2e^{u_2}}{\int_\Si h_2e^{u_2}dV_g}-1\rr-\rho_1\lr\fr{h_1e^{u_1}}{\int_\Si h_1e^{u_1}dV_g}-1\rr-4\pi\sum_{j=1}^J\a_{2j}\lr\d_{p_j}-1\rr\ea\ry,$$
where $\dis{h_i}$ are smooth positive functions, $\dis{\rho_i}$ are positive real parameters, $\dis{p_j}$ are given points on $\dis{\Si}$ and $\dis{\a_{ij}}$ are numbers greater than $\dis{-1}$.\\
We give existence and multiplicity results, using variational and Morse-theoretical methods. It is the first existence result when some of the $\dis{\a_{ij}}$'s are allowed to be negative.
\eab\

\sec{Introduction}\

Let $\dis{\Si}$ be a compact surface without boundary and $\dis{g}$ a Riemannian metric on $\dis{\Si}$. The $\dis{SU(N+1)}$ Toda system is the following system of elliptic PDEs:
\beq\lab{todaN}
-\D u_i=\sum_{j=1}^Na_{ij}\rho_j(h_je^{u_j}-1),\q\q\q i=1,\ds,N
\eeq
where $\dis{\D=\D_g}$ is the Laplace-Beltrami operator, $\dis{\rho_i}$ are positive real parameters, $\dis{h_i}$ are smooth positive functions and $\dis{A=(a_{ij})_{ij}}$ is the Cartan matrix of $\dis{SU(N+1)}$
$$\lr\ba{ccccc}2&-1&0&\cds&0\\-1&2&\dds&\dds&\vds\\0&\dds&\dds&\dds&0\\\vds&\dds&\dds&2&-1\\0&\cds&0&-1&2\ea\rr.$$
The Toda system has been widely studied due to its great importance in both geometry and mathematical physics: in geometry, it arises in the description of holomorphic curves in $\dis{\C\P^N}$ (see e.g. \ci{bw,cal,cw}), whereas in mathematical physics it is a model for non-abelian Chern-Simons vortices theory (see \ci{dunne,tar08,yang}).\\
It is not restrictive to suppose the total area $\dis{|\Si|}$ of $\dis{\Si}$ to be equal to $\dis{1}$; therefore, integrating $\dis{\eqr{todaN}}$ on $\dis{\Si}$, we deduce that any solution verifies
$$\int_\Si h_ie^{u_i}dV_g=1\q\q\q\fa\,i=1,\ds,N,$$
hence the system $\dis{\eqr{todaN}}$ is equivalent to
$$-\D u_i=\sum_{j=1}^Na_{ij}\rho_j\lr\fr{h_je^{u_j}}{\int_\Si h_je^{u_j}dV_g}-1\rr,\q\q\q i=1,\ds,N$$
which has the advantage of being invariant by addition of constants.\\

A variant of this system is given by adding in the right-hand side a linear combination of Dirac deltas centered at points of $\dis{\Si}$.\\
This variant has still applications in mathematical physics and geometry. In the former, it arises in gauged self-dual Schr\"odinger equations (see \ci{yang}), where the supports of Dirac deltas represent the \qm{vortices} of the wave function, that is the points where it vanishes. In the latter, it is related to the study of holomorphic curves with ramifications: here, the ramificated points are the centers of the Dirac deltas and the ramification index is given by the coefficient multiplying the delta.\\
In particular, in this paper we will consider the $\dis{SU(3)}$ Toda system with singularities
\beq\lab{toda2}
\lb\bl-\D u_1=2\rho_1\lr\fr{h_1e^{u_1}}{\int_\Si h_1e^{u_1}dV_g}-1\rr-\rho_2\lr\fr{h_2e^{u_2}}{\int_\Si h_2e^{u_2}dV_g}-1\rr-4\pi\sum_{j=1}^J\a_{1j}\lr\d_{p_j}-1\rr\\-\D u_2=2\rho_2\lr\fr{h_2e^{u_2}}{\int_\Si h_2e^{u_2}dV_g}-1\rr-\rho_1\lr\fr{h_1e^{u_1}}{\int_\Si h_1e^{u_1}dV_g}-1\rr-4\pi\sum_{j=1}^J\a_{2j}\lr\d_{p_j}-1\rr\ea\ry.
\eeq\\

To better describe the main properties of this system, let us perform a change of variables. Consider the Green function $\dis{G_p}$ of the Laplace operator centered at a point $\dis{p\in\Si}$, that is the solution of
$$\lb\bl-\D G_p=\d_p-1\\\int_\Si G_pdV_g=0\ea\ry,$$
and apply the change of variables
\beq\lab{cv}
u_i\to u_i+4\pi\sum_{j=1}^J\a_{ij}G_{p_j}.
\eeq
Then problem $\dis{\eqr{toda2}}$ transforms into the following:
\beq\lab{toda}
\lb\bl-\D u_1=2\rho_1\lr\fr{\wt h_1e^{u_1}}{\int_\Si\wt h_1e^{u_1}dV_g}-1\rr-\rho_2\lr\fr{\wt h_2e^{u_2}}{\int_\Si\wt h_2e^{u_2}dV_g}-1\rr\\-\D u_2=2\rho_2\lr\fr{\wt h_2e^{u_2}}{\int_\Si\wt h_2e^{u_2}dV_g}-1\rr-\rho_1\lr\fr{\wt h_1e^{u_1}}{\int_\Si\wt h_1e^{u_1}dV_g}-1\rr\ea\ry.
\eeq
Here, the functions $\dis{\wt h_i}$ are defined by
$$\wt h_i=h_ie^{-4\pi\sum_{j=1}^J\a_{ij}G_{p_j}},\q i=1,2,$$
therefore they verify
$$\lb\bll0<\wt h_i\in C^\ity\lr\Si\bs\Cup_{j=1}^Jp_j\rr\\\wt h_i\sim d(\cd,p_j)^{2\a_{ij}}&\tx{near }p_j\ea\ry,$$
hence $\dis{\wt h_i}$ has a singularity at $\dis{p_j}$ if $\dis{\a_{ij}<0}$ and it has a zero at $\dis{p_j}$ if $\dis{\a_{ij}>0}$.\\

Problem $\dis{\eqr{toda}}$ has a variational formulation, that is its solutions are the critical points of the Euler-Lagrange functional defined by
\beq\lab{jrho}
J_\rho(u)=\int_\Si Q(u)dV_g-\sum_{i=1}^2\rho_i\lr\log\int_\Si\wt h_ie^{u_i-\int_\Si u_idV_g}dV_g\rr,
\eeq
with $\dis{\rho=(\rho_1,\rho_2)}$, $\dis{u=(u_1,u_2)}$ and $\dis{Q(u)}$ is given by
$$Q(u)=\fr{1}3\lr|\n u_1|^2+\n u_1\cd\n u_2+|\n u_2|^2\rr,$$
where $\dis{\n=\n_g}$ is the gradient given by the metric $\dis{g}$ and $\dis{\cd}$ denotes the Riemannian scalar product.\\
A first tool to study the structure of the functional $\dis{J_\rho}$ is given by the following Moser-Trudinger inequality, which was proved in \ci{batmal} (and, for the regular case, in \ci{jw}):
\beq\lab{mttoda}
\sum_{i=1}^2\min\lb1,1+\min_j\a_{ij}\rb\lr\log\int_\Si\wt h_ie^{u_i-\int_\Si u_idV_g}dV_g\rr\le\fr{1}{4\pi}\int_\Si Q(u)dV_g+C.
\eeq
This inequality implies boundedness from below whenever $\dis{\rho_i\le4\pi\min\lb1,1+\min_j\a_{ij}\rb}$ for both $\dis{i=1,2}$ and coercivity (up to addition of constants) if both $\dis{\rho_i}$'s are strictly smaller, therefore in this case $\dis{\eqr{toda}}$ has a minimizing solution. On the other hand, in the same papers it is shown, through suitable test functions, that $\dis{J_\rho}$ is unbounded from below for greater values of $\dis{\rho_i}$, so one can no longer find critical points through minimization techniques.\\

The first main result of this paper is about existence of solutions on surfaces with positive genus and arbitrarily signed vortices.\\
Before stating it, let us apply a change of notation about the singular points and their singularities. Since we will suppose to have $\dis{\max\{\a_{1j},\a_{2j}\}\ge0}$, we can divide the singular points into three categories, depending on whether the first component, the second component or none of them has a negative singularities on it. We also consider, alongside the multi-indices $\dis{\a_1,\a_2}$, two sub-indices $\dis{\wt\a_1}$, $\dis{\wt\a_2}$ which take account only of the negative $\dis{\a_{ij}}$'s, and we order them in such a way that they are not decreasing.\\
Precisely, we write
$$\{p_1,\ds,p_J\}=\lb p_{01},\ds,p_{0L_0},p_{11},\ds,p_{1L_1},p_{21},\ds,p_{2L_2}\rb$$
with $\dis{p_j=p_{il}}$ for some $\dis{i=1,2}$, $\dis{l=1,\ds,L_i}$ if and only if $\dis{\wt\a_{il}:=\a_{ij}<0}$ and $\dis{\wt\a_{i1}\le\ds\le\wt\a_{iL_i}}$.\\

\bth\lab{ex}$\dis{}$\\
Suppose $\dis{\Si}$ has positive genus and $\dis{\max\{\a_{1j},\a_{2j}\}\ge0}$ for any $\dis{j=1,\ds,J}$.
Then, there exists a closed set $\dis{\L\sub\R_+^2}$ with zero Lebesgue measure such that for any $\dis{\rho=(\rho_1,\rho_2)\nin\L}$ which satisfies
\beq\lab{rhoi}
4\pi\lr K_i+\sum_{l\in\mc I_i}(1+\wt\a_{il})\rr<\rho_i<4\pi\lr K_i+\sum_{l\in\mc I_i\cup\{1\}}(1+\wt\a_{il}),\rr\q i=1,2
\eeq
for some $\dis{K_i\in\N}$ and $\dis{\mc I_i\sub\{2,\ds,L_i\}}$ the problem $\dis{\eqr{toda}}$ admits at least one solution.
\eth\

In the last section we will give some examples to clarify the meaning of condition $\dis{\eqr{rhoi}}$.\\

Theorem $\dis{\ref{ex}}$ is, up to our knowledge, the first existence result for the singular Toda system with arbitrarily signed vortices. A recent paper \ci{bjmr} gives a general existence result for $\dis{\eqr{toda2}}$ when $\dis{\a_{ij}\ge0}$ and $\dis{g(\Si)>0}$, and there are some other partial existence results for the regular case, i.e. with upper bounds on one or both of the $\dis{\rho_i}$'s (\ci{mn,mr13}).\\
Removing the hypothesis of non-negativity of the vortices reduces the similarities with the regular case - for instance the best constant in the Moser-Trudinger inequality $\dis{\eqr{mttoda}}$ is no longer the same - thus increasing the difficulty of the problem.\\
This issue also arose clearly in the study of the scalar counterpart of $\dis{\eqr{toda2}}$, that is the singular Liouville equation
$$-\D u=2\rho\lr\fr{he^u}{\int_\Si he^udV_g}-1\rr-4\pi\sum_{j=1}^J\a_j\lr\d_{p_j}-1\rr.$$
which, by a trick similar to $\dis{\eqr{cv}}$, is equivalent to
\beq\lab{liou}
-\D u=2\rho\lr\fr{\wt he^u}{\int_\Si\wt he^udV_g}-1\rr.
\eeq\\

Equation $\dis{\eqr{liou}}$ is also very important in both geometry and mathematical physics: it arises in the problem of prescribed Gauss curvature on surfaces with conical singularities and it appears in some models in Chern-Simons theory. It has been widely studied in literature, with many results on existence, compactness of solutions etc., which have been reviewed in \ci{mal10,tar10}.\\
In the scalar case, general existence results have been found in the case of positive genus and non-negative vortices, both through variational and Morse-theoretical methods \ci{bdm} and through the computation of the Leray-Schauder degree \ci{cl13}. On the other hand, when the coefficients $\dis{\a_j}$ may attain negative values, the existence or non-existence of solutions depends on $\dis{\rho}$ and on the $\dis{\a_j}$'s (\ci{car,carmal}), as well as for the general case of the sphere (\ci{barmal,mr11}).\\

In this paper, we also give a generic multiplicity result for the problem $\dis{\eqr{toda}}$, using Morse theory. Basically, the more are the couples $\dis{(K,\mc I)}$ which satisfy $\dis{\eqr{rhoi}}$, the higher the number of solution is.\\

\bth\lab{mult}$\dis{}$\\
Suppose the hypotheses of Theorem $\dis{\ref{ex}}$ hold, and suppose that for $\dis{i=1,2}$ there exist $\dis{H_i,K_{i1},\ds,K_{iH_i}\in\N}$ and $\dis{\mc I_{i1},\ds,\mc I_{iH_i}\sub\{2,\ds,L_i\}}$ such that any $\dis{h=1,\ds,H_i}$ verifies
$$4\pi\lr K_{ih}+\sum_{l\in\mc I_{ih}}(1+\wt\a_{il})\rr<\rho_i<4\pi\min\lb K_{ih}+\sum_{l\in\mc I_{ih}\cup\{1\}}(1+\wt\a_{il}),K_{ih}+1+\sum_{l\in\mc I_{ih}\bs\lb\max\mc I_{ih}\rb}(1+\wt\a_{il})\rb.$$
Then, there exists a dense open set of $\dis{D\sub\mc M^2(\Si)\x L^\ity(\Si)^2}$ such that if $\dis{(g,h_1,h_2)\in D}$, then the problem $\dis{\eqr{toda}}$ has at least
$$\sum_{h_1,h_2}\bin{K_{1h_1}+|\mc I_{1h_1}|+\ls\fr{-\chi(\Si)}2\rs}{|\mc I_{1h_1}|+\ls\fr{-\chi(\Si)}2\rs}\bin{K_{2h_2}+|\mc I_{2h_2}|+\ls\fr{-\chi(\Si)}2\rs}{|\mc I_{2h_2}|+\ls\fr{-\chi(\Si)}2\rs}$$
solutions,
where $\dis{\mc M^2(\Si)}$ is the space of Riemannian metrics on $\dis{\Si}$ endowed with the $\dis{C^2}$ topology and square brackets denote the integer part of a real number.
\eth\

We stress that, up to our knowledge, there is no previous multiplicity result for the Toda system, even in the regular case. For the Liouville equation multiplicity results have been obtained using both Morse theory (\ci{bdm,dem}) and topological degree (\ci{cl03,cl13}).\\
In particular, the multiplicity result for the case of non-negative singularities has a quite simpler form, which is summarized in the following corollary:\\

\bco$\dis{}$\\
Suppose $\dis{g(\Si)>0}$, $\dis{\a_{ij}\ge0}$ for all $\dis{i,j}$ and
$$\rho\in(4K_1\pi,4(K_1+1)\pi)\x(4K_2\pi,4(K_2+1)\pi)\bs\L.$$
Then, for a generic choice of the data (as in Theorem $\dis{\ref{mult}}$) the problem $\dis{\eqr{toda}}$ admits at least
$$\bin{K_1+\ls\fr{-\chi(\Si)}2\rs}{\ls\fr{-\chi(\Si)}2\rs}\bin{K_2+\ls\fr{-\chi(\Si)}2\rs}{\ls\fr{-\chi(\Si)}2\rs}$$
solutions.
\eco\

The same arguments of Theorems $\dis{\ref{ex}}$ and $\dis{\ref{mult}}$ can also be applied in a couple of other cases.\\
First of all, considering again surfaces with positive genus, we can remove the hypotheses $\dis{\max\{\a_{1j},\a_{2j}\}\ge0}$ if we suppose at least one of the parameter $\dis{\rho_i}$ to be small enough, that is for instance if $\dis{\rho_2<4\pi(1+\a_{2j})}$ for all $\dis{j}$'s such that both $\dis{\a_{1j}}$ and $\dis{\a_{2j}}$ are negative (hence, in particular, if it satisfies the coercivity condition for $\dis{\eqr{mttoda}}$).\\
Moreover, on surfaces of arbitrary genus, we can obtain a similar result concerning both existence and multiplicity of solutions if we assume both $\dis{\rho_1}$ and $\dis{\rho_2}$ to satisfy the upper bound stated before.\\
In both cases, we will again consider points $\dis{p_{0l},p_{1l},p_{2l}}$ and sub-indices $\dis{\wt\a_1,\wt\a_2}$ as before, though considering only the negative $\dis{\a_{ij}}$'s which can be attained in the restricted range of $\dis{\rho_i}$, that is $\dis{\a_{2j}<\wh\a_2}$ in the former case and $\dis{\a_{ij}<\wh\a_i}$ for both $\dis{i}$ in the latter case, where
\beq\lab{aiji}
\wh\a_i:=\inf\{\a_{ij}:\max\{\a_{1j},\a_{2j}\}<0\}.
\eeq

\bth\lab{exmult}$\dis{}$\\
Suppose $\dis{g(\Si)>0}$ and that $\dis{\rho\in\R^2_+\bs\L}$ verifies $\dis{\rho_2<4\pi(1+\wh\a_i)}$, with $\dis{\wh\a_i}$ as in $\dis{\eqr{aiji}}$, and
$$4\pi\lr K+\sum_{l\in\mc I_1}(1+\wt\a_{1l})\rr<\rho_1<4\pi\lr K+\sum_{l\in\mc I_1\cup\{1\}}(1+\wt\a_{1l})\rr$$
\beq\lab{rho2}
4\pi\sum_{l\in\mc I_2}(1+\wt\a_{2l})<\rho_2<4\pi\sum_{l\in\mc I_2\cup\{1\}}(1+\wt\a_{2l})
\eeq
for some $\dis{K\in\N}$ and $\dis{\mc I_i\sub\{1,\ds,L_i\}}$. Then, the problem $\dis{\eqr{toda}}$ admits at least a solution.\\
If moreover the condition $\dis{\eqr{rho2}}$ is satisfied by $\dis{\mc I_{21},\ds,\mc I_{2H_2}}$ and there exist $\dis{H_1,K_1,\ds,K_{H_1}\in\N}$ and $\dis{\mc I_{11},\ds,\mc I_{1H_1}\sub\{2,\ds,L_1\}}$ satisfying, for any $\dis{h=1,\ds,H_1}$,
$$4\pi\lr K_h+\sum_{l\in\mc I_{1h}}(1+\wt\a_{1l})\rr<\rho_1<4\pi\min\lb K_h+\sum_{l\in\mc I_{1h}\cup\{1\}}(1+\wt\a_{1l}),K_h+1+\sum_{l\in\mc I_{1h}\bs\lb\max\mc I_{1h}\rb}(1+\wt\a_{1l})\rb,$$
then a generic choice of data yields at least
$$H_2\sum_{h}\bin{K_h+|\mc I_{1h}|+\ls\fr{-\chi(\Si)}2\rs}{|\mc I_{1h}|+\ls\fr{-\chi(\Si)}2\rs}$$
solutions.
\eth

\bth\lab{sph}$\dis{}$\\
Suppose $\dis{\rho\in\R_+^2\bs\L}$ verifies $\dis{\rho_i<4\pi\min\{1,1+\wh\a_i\}}$ for both $\dis{i=1,2}$ and
$$4\pi\sum_{l\in\mc I_i}(1+\wt\a_{il})<\rho_i<4\pi\sum_{l\in\mc I_i\cup\{1\}}(1+\wt\a_{il}),\q i=1,2$$
for some $\dis{\mc I_i\sub\{2,\ds,L_i\}}$. Then, the problem $\dis{\eqr{toda}}$ admits at least a solution.
Moreover, if the above condition is verified for $\dis{\mc I_{11},\ds,\mc I_{1H_1}}$ and $\dis{\mc I_{21},\ds,\mc I_{2H_2}}$, then a generic choice of initial data yields at least $\dis{H_1H_2}$ solutions.
\eth\

The set $\dis{\L}$ in the statement of Theorem $\dis{\ref{ex}}$ can be explicitly written as a union of straight lines and points in dependence of the $\dis{\a_{ij}}$ (see the next section) and it arises in the study of compactness and blow-up of solutions of $\dis{\eqr{toda}}$. This has been one of the major difficulties in attacking both $\dis{\eqr{toda}}$ and $\dis{\eqr{liou}}$.\\
In the scalar case, quantization results have been given (see \ci{bremer,li,ls} for the regular and \ci{barmon,bt2,bt1} for the singular case): a sequence $\dis{\{u_n\}_{n\in\N}}$ of solutions of $\dis{\eqr{liou}}$ which blows up at a regular point $\dis{p\nin\{p_1,\ds,p_J\}}$ satisfies
$$\lim_{r\to0}\lim_{n\to+\ity}\int_{B_r(p)}\wt he^{u_n}dV_g=4\pi,$$
whereas if it blows up at $\dis{p_j}$ it verifies
$$\lim_{r\to0}\lim_{n\to+\ity}\int_{B_r(p_j)}\wt he^{u_n}dV_g=4\pi(1+\a_j).$$
Much is also known about the blow-up behavior of the regular Toda system.\\
In \ci{jlw} it was proved that there are basically three different blow-up scenarios. The first occurs when only one component $\dis{u_i}$ is blowing up: in this case the quantization values for the two components are $\dis{(4\pi,0)}$ or $\dis{(0,4\pi)}$. If both components blow up at different rates of concentration, the quantization values are $\dis{(8\pi,4\pi)}$ (respectively $\dis{(4\pi,8\pi)}$), whereas if the two components blow up at the same rate it is $\dis{(8\pi,8\pi)}$. In $\ci{egp,mpw}$ it was shown that all these cases are actually possible.\\
In the presence of singularities, the expected corresponding quantization values at a singular point $\dis{p_j}$ would be respectively
$$(4\pi(1+\a_{1j}),0),\q(0,4\pi(1+\a_{2j})),\q(4\pi(1+\a_{1j}),4\pi(2+\a_{1j}+\a_{2j}))$$
$$(4\pi(2+\a_{1j}+\a_{2j}),4\pi(1+\a_{2j})),\q(4\pi(2+\a_{1j}+\a_{2j}),4\pi(2+\a_{1j}+\a_{2j})).$$
Blow-up phenomena for $\dis{\eqr{toda}}$ have been investigated in \ci{lwz}: the authors showed that only a finite number of blow-up values are allowed, including the five above (see the next section for details). However, it is an open problem whether these values indeed occur or whether they can be excluded as well.\\

Let us see now the role played by the study of sub-levels of the energy functional for the existence of solutions for $\dis{\eqr{toda}}$.\\
Concerning the scalar Liouville equation $\dis{\eqr{liou}}$, the Euler functional
\beq\lab{irho}
I_\rho(u)=\fr{1}2\int_\Si|\n u|^2dV_g-2\rho\log\int_\Si\wt he^{u-\int_\Si udV_g}dV_g
\eeq
is bounded from below if and only if $\dis{\rho\le4\pi\min\lb1,1+\min_j\a_j\rb}$ and it is coercive if and only if $\dis{\rho}$ is smaller, as follows from the inequalities in \ci{fon,mos,tro}.\\
To study variationally the problem for higher values of the parameter $\dis{\rho}$, a first clue was given by Chen and Li \ci{cl01}, who showed that $\dis{I_\rho}$ is bounded from below under the assumption of some spreading of the ($\dis{L^1}$-normalized) function $\dis{\wt h e^u}$. Improving their result, Djadli \ci{dja} and Malchiodi \ci{mal08} gave a general existence result for the regular case of $\dis{\eqr{liou}}$ by showing, when $\dis{\rho\in(4K\pi,4(K+1)\pi)}$, a homotopy equivalence between the low sub-levels of $\dis{I_\rho}$ and the non-contractible set of formal barycenters on $\dis{\Si}$
\beq\lab{bar}
\Si^K:=\lb\sum_{k=1}^Kt_k\d_{x_k};\,x_k\in\Si,\,t_k\ge0,\,\sum_{k=1}^Kt_k=1\rb.
\eeq
An extension of this was later given by Carlotto and Malchiodi \ci{carmal} who considered the Liouville equation with non-positive singularities. The authors extended the notion of formal barycenters by defining a sort of weighted ones on $\dis{\Si}$: since the constant in Moser-Trudinger inequality worsens near the singularities, they redefined the set $\dis{\eqr{bar}}$ in such a way that the points $\dis{p_j}$ are somehow \qm{lighter} than regular ones, proportionally to the respective coefficient $\dis{\a_j}$. They defined the weight of a finite set of $\dis{\Si}$ with respect to the multi-index $\dis{\a=(\a_1,\ds,\a_J)}$ as
\beq\lab{weight}
\mc J=\{q_1,\ds,q_K,p_{j_1},\ds,p_{j_L}\}\q\To\q\chi_\a(\mc J)=K+\sum_{l=1}^L(1+\a_{j_l})
\eeq
\beq\lab{wbar}
\Si_{\rho,\a}=\lb\sum_{x_k\in\mc J}t_k\d_{x_k};\,x_k\in\Si,\,t_k\ge0,\,\sum_{x_k\in\mc J}t_k=1,\,4\pi\chi_\a(\mc J)<\rho\rb.
\eeq
They showed that the homology groups of $\dis{\Si_{\rho,\a}}$ are mapped invectively into the ones of $\dis{I_\rho}$'s sub-levels. However, the topological structure of the weighted barycenters depends heavily on the parameters $\dis{\rho}$ and $\dis{\a_j}$'s and it can be much more complicated than in the regular case; this is discussed in \cite{car}.\\

Concerning the regular Toda system, the argument in \ci{dja,mal08} was adapted in \ci{mn} to the case of $\dis{\rho_1<4\pi}$ and $\dis{\rho_2\in(4K\pi,4(K+1)\pi)}$, since the second component has the same concentration behavior which occurs in the scalar case.\\
If instead both parameters are supercritical, both components can concentrate, thus interacting in a definitely non-trivial way, as discussed in \ci{mr13} for $\dis{\rho\in(4\pi,8\pi)^2}$.\\
Another difficulty which might arise in the singular case is the concentration around positively-signed vertices. The presence of singularities affects significantly the bubbling behavior but at the same time it would make no sense to assign them a different weight from the regular points, since they make no difference for what concerns the constants in the Moser-Trudinger inequality.\\
To overcome these difficulties, we adapt a topological argument from \ci{bdm,bjmr}. Since we suppose $\dis{\Si}$ neither being homeomorphic to $\dis{\S^2}$ nor to $\dis{\R\P^2}$, we can take two bouquets of circles $\dis{\g_1,\g_2}$ (that is, two collection of circles glued each around a single point) such that $\dis{\Si}$ can retract on each of them through continuous maps $\dis{\Pi_1,\Pi_2}$. For our purpose, we choose $\dis{\g_1}$ containing all the points $\dis{p_{11},\ds,p_{1L_1}}$ (using the same notation as in Theorem $\dis{\ref{ex}}$) and none of the other singular points, and in the same way we choose $\dis{\g_2}$ containing, among the singular points, all and only the $\dis{p_{2l}}$'s.\\
If we are under the hypotheses of Theorem $\dis{\ref{exmult}}$, by the different way we defined the points $\dis{p_{il}}$, there might still be concentration of both components around the same negative singularity, but this is actually excluded by assuming $\dis{\rho_2<4\pi(1+\a_{2j})}$.\\
The same difficulties can be similarly avoided, in Theorem $\dis{\ref{sph}}$, even in the zero-genus case. In fact, thanks to the hypothesis $\dis{\rho_1,\rho_2<4\pi}$, as in Theorem $\dis{\ref{sph}}$, each component can concentrate only around negative singularities, so the interaction does not occur because we also assumed $\dis{\rho_i<4\pi(1+\wh\a_i)}$.\\
Putting together arguments from \ci{bjmr,carmal,mn,mr13} we show that if $\dis{J_\rho(u)}$ is sufficiently low, then for one or both $\dis{i=1,2}$, $\dis{\wt h_ie^{u_i}}$ is arbitrarily close (in some sense which will be better specified in the next sections) to the set $\dis{\Si_{\rho_i,\wt\a_i}}$, defined as $\dis{\eqr{wbar}}$ with multi-index $\dis{\wt\a_i=\lr\wt\a_{i1},\ds,\wt\a_{iL_i}\rr}$. Therefore, it is possible to map continuously low sub-levels of $\dis{J_\rho}$ on one or both the $\dis{\Si_{\rho_i,\wt\a_i}}$'s and, through the retractions $\dis{\Pi_i}$'s, on $\dis{(\g_i)_{\rho_i,\wt\a_i}}$.\\
To express the fact that only one or both mappings can be built, as in \ci{bjmr} we introduce the topological join $\dis{(\g_1)_{\rho_1,\wt\a_1}\st(\g_2)_{\rho_2,\wt\a_2}}$: the topological join $\dis{X\st Y}$ of two sets $\dis{X}$ and $\dis{Y}$ is basically the product $\dis{X\x Y\x[0,1]}$ with the endpoints $\dis{X\x Y\x\{0\}}$ and $\dis{X\x Y\x\{1\}}$ collapsed respectively on $\dis{X\x\{0\}}$ and $\dis{Y\x\{1\}}$.\\
Therefore, we are able to define a projection $\dis{\Psi}$ from low sub-levels of $\dis{J_\rho}$ to the join $\dis{(\g_1)_{\rho_1,\wt\a_1}\st(\g_2)_{\rho_2,\wt\a_2}}$. On the other hand, we can also build a map $\dis{\Phi}$ from the join of the weighted barycenter sets to arbitrarily low sub-levels of $\dis{J_\rho}$ through suitable test functions.\\
The composition $\dis{\Psi\c\Phi}$ is homotopically equivalent to the identity on the join of the barycenters, so the homology groups of $\dis{J_\rho}$'s sub-levels contain a copy of the ones of the join. Finally, since we can express with a simple formula the homology of $\dis{A\st B}$ in terms of the homology of $\dis{A}$ and $\dis{B}$, we deduce non-contractibility of low sub-levels for suitable values of $\dis{\rho_i}$ and $\dis{\a_{ij}}$.\\
To conclude the proof of the existence result, we would need a Palais-Smale-like compactness condition. Palais-Smale condition is not known to hold for $\dis{J_\rho}$ but for a dense set of the parameters $\dis{\rho}$ bounded Palais-Smale sequences exist, as follows from \ci{str}. Therefore, we obtain existence of solution for $\dis{\rho}$ belonging to a dense set; the compactness of solutions (which follows from assuming $\dis{\rho\nin\L}$) permits to extend the existence result to any admissible $\dis{\rho}$.\\

To get a multiplicity result, we use the weak Morse inequalities, which relate the number of critical points to the homology groups.\\
In particular, we get that the total number of solutions of $\dis{\eqr{toda}}$, that is the critical points of $\dis{J_\rho}$, is greater or equal to the Betti numbers of $\dis{J_\rho}$'s low sub-levels.\\
By the above analysis, the latter will be greater or equal to the ones of $\dis{(\g_1)_{\rho_1,\wt\a_1}\st(\g_2)_{\rho_2,\wt\a_2}}$, which can be estimated from below using the Mayer-Vietoris exact sequence.\\
The fact that $\dis{J_\rho}$ is a Morse function for a generic choice of $\dis{g}$, $\dis{h_i}$, $\dis{p_j}$ and $\dis{\a_{ij}}$ follows by arguing as in \ci{bdm,dem}.\\

In Section $\dis{2}$ we introduce some notation and some preliminary results which will be used later on. In Section $\dis{3}$ we build the map $\dis{\Phi}$ from the join of the weighted barycenter sets to low sub-levels of $\dis{J_\rho}$. Section $\dis{4}$ is devoted to the study of the variational structure of $\dis{J_\rho}$ and to obtaining improved Moser-Trudinger inequalities; the latter results will be needed to construct the map $\dis{\Psi}$ from sub-levels of $\dis{J_\rho}$ to $\dis{(\g_1)_{\rho_1,\wt\a_1}\st(\g_2)_{\rho_2,\wt\a_2}}$ and to prove that the composition with $\dis{\Phi}$ is homotopically equivalent to the identity, which is done in Section $\dis{5}$. In Section $\dis{6}$ we study the topology and the homology of the set of weighted barycenters. Finally, in Section $\dis{7}$ we see some examples and we put together the result obtained to prove the existence and multiplicity results.\\

\sec{Notation and preliminaries}\

In this section we will provide some notation and some known preliminary results which we will need later.\\
The indicator function of a set $\dis{\O\sub\Si}$ will be denoted as
$$\mf 1_\O(x)=\lb\bll1&\tx{if }x\in\O\\0&\tx{if }x\nin\O\ea\ry$$
Given two points $\dis{x,y\in\Si}$, we will denote the metric distance between them on $\dis{\Si}$ as $\dis{d(x,y)}$. In the same way, for any two subsets $\dis{\O,\O'\sub\Si}$ we will denote:
$$d(x,\O):=\inf\{d(x,y):y\in\O\},\q\q\q d(\O,\O'):=\inf\lb d(x,y):x\in\O,y\in\O'\rb.$$
We will denote as $\dis{D_\Si}$ the diameter of $\dis{\Si}$
$$D_\Si:=\sup\{d(x,y):\,x,y\in\Si\}.$$
The symbol $\dis{B_r(p)}$ will stand for the open metric ball centered at $\dis{p}$ and having radius $\dis{r}$. We will similarly use the notation $\dis{B_r(\O)}$ for a subset $\dis{\O\sub\Si}$:
$$B_r(\O):=\{x\in\Si:d(x,\O)<r\}.$$
Given a function $\dis{u\in L^1(\Si)}$ and a measurable set $\dis{\O\sub\Si}$, we denote the average of $\dis{u}$ on $\dis{\O}$ as
$$\int_\O udV_g=\fr{1}{|\O|}\int_\O udV_g$$
The symbol $\dis{\ov u}$ will stand for the average of $\dis{u}$ on $\dis{\Si}$; since we assume $\dis{|\Si|=1}$, we can write
$$\ov u=\int_\Si udV_g=\fin_\Si udV_g.$$
We will denote the subset of functions in $\dis{H^1(\Si)}$ having null average as
$$\ov H^1(\Si):=\lb u\in H^1(\Si):\,\ov u=0\rb.$$
Notice that, since the functional $\dis{J_\rho}$ defined in $\dis{\eqr{jrho}}$ is invariant by addition of constants, it will not be restrictive to study it on $\dis{\ov H^1(\Si)^2}$ rather than on $\dis{H^1(\Si)^2}$.\\
The sub-levels of $\dis{J_\rho}$, which, as anticipated, will play and essential role throughout the whole paper, will be denoted as
$$J_\rho^a=\lb u\in H^1(\Si)^2:\,J_\rho(u)\le a\rb.$$
For a continuous map $\dis{f:\Si\to\Si}$ and a measure $\dis{\mu\in\mc M(\Si)}$, we define the push-forward of $\dis{\mu}$ with respect to $\dis{f}$ the measure defined by
$$f_*\mu(B)=\mu\lr f^{-1}(B)\rr.$$
If $\dis{\mu}$ has finite support, its push-forward has a particularly simple form:
$$\mu=\sum_{k=1}^Kt_k\d_{x_k}\q\q\q\To\q\q\q f_*\mu=\sum_{k=1}^Kt_k\d_{f(x_k)}.$$
We will use the symbol $\dis{X\seq Y}$ to mean that two topological spaces $\dis{X}$ and $\dis{Y}$ are homotopically equivalent. We will consider the composition of two homotopy equivalence $\dis{H_1:X\x[0,1]\to Y}$, $\dis{H_2:Y\x[0,1]\to Z}$, that is the map $\dis{H_2*H_1:X\x[0,1]\to Z}$ defined by
$$H_2*H_1:(x,s)\to\lb\bll H_1(x,2s)&\tx{if }s\le\fr{1}2\\H_2(x,2s-1)&\tx{if }s>\fr{1}2\ea\ry.$$
The identity map on $\dis{X}$ will be denoted by $\dis{\mr{Id}_X}$.
Given a topological space $\dis{X}$ we will denote its $\dis{q^\mr{th}}$ homology group with coefficients in $\dis{\Z}$ as $\dis{H_q(X)}$. Isomorphisms between homology group will be denoted just by the equality sign. We will denote as $\dis{\wt H_q(X)}$ the $\dis{q}$-th reduced homology group with coefficients in $\dis{\Z}$, that is
$$H_0(X)=\wt H_0(X)\pl\Z,\q\q\q H_q(X)=\wt H_q(X)\q\tx{if }q\ne0$$
The symbol $\dis{\b_q(X)}$ will stand for the $\dis{q}$-th Betti number of $\dis{X}$, that is $\dis{\b_q(X)=\mr{rank}(H_q(X))}$. As before, $\dis{\wt\b_q(X)}$ will stand for the dimension of $\dis{\wt H_q(X)}$, therefore it will coincide with the usual definition of Betti number with the exception of $\dis{\wt\b_0(X)=\b_0(X)-1}$. For a subspace $\dis{Y\sub X}$, the $\dis{q}$-th relative homology group (with coefficients in $\dis{\Z}$) will be denoted by $\dis{H_q(X,Y)}$ and $\dis{\b_q(X,Y)}$ will be the relative Betti numbers.\\
If $\dis{J_\rho}$ is a Morse function, we will denote as $\dis{\mc C_q(a,b)}$ the number of critical points $\dis{u}$ of $\dis{J_\rho}$ with Morse index $\dis{q}$ satisfying $\dis{a\le J_\rho(u)\le b}$. The total number of critical points of index $\dis{q}$ will be denoted as $\dis{\mc C_q}$; in other words, $\dis{\mc C_q:=\mc C_q(+\ity,-\ity)}$.\\

Throughout all the paper we will denote by $\dis{C}$ large constants which can vary among different lines or formulas. To stress the dependence of $\dis{C}$ on some parameter we will add subscripts such as $\dis{C_\a}$ and so on.\\
We will use the symbol $\dis{o_\a(1)}$ to denote quantities which tend to $\dis{0}$ as $\dis{\a}$ goes to $\dis{0}$ or to $\dis{+\ity}$ and we will similarly write $\dis{O_\a(1)}$ for bounded quantities. The subscript will be omitted when it is evident from the context.\\
In a similar way, we will use the symbol $\dis{f_\a\sim_\a g_\a}$, or simply $\dis{\sim}$, to express that the ratio between $\dis{f_\a}$ and $\dis{g_\a}$ is bounded by a positive constant both from above and below when $\dis{\a}$ goes to $\dis{0}$ or $\dis{+\ity}$. In other words, this means that $\dis{\log\fr{f_\a}{g_\a}=O_\a(1)}$.\\

Now we recall the Moser-Trudinger inequality for the Liouville equation and the Toda system and their immediate corollaries.\\

\bth(\ci{fon}, Theorem $\dis{1.7}$; \ci{mos}, Theorem $\dis{2}$; \ci{tro}, Corollary $\dis{9}$)$\dis{}$\\
For any $\dis{u\in H^1(\Si)}$ it holds
\beq\lab{mtliou}
\log\int_\Si\wt he^{u-\ov u}dV_g\le\fr{1}{16\pi\min\{1,1+\min_j\a_j\}}\int_\Si|\n u|^2dV_g+C.
\eeq
In other words, the functional $\dis{I_\rho}$ defined in $\dis{\eqr{irho}}$ is bounded from below if and only if $\dis{\rho\le4\pi\min\lb1,1+\min_j\a_j\rb}$.
\eth\

\bco$\dis{}$\\
The functional $\dis{I_\rho}$ is coercive on $\dis{\ov H^1(\Si)}$ if and only if  $\dis{\rho<4\pi\min\lb1,1+\min_j\a_j\rb}$. If this occurs, then it has a global  minimizer $\dis{u}$ which solves $\dis{\eqr{liou}}$.
\eco\

\bth(\ci{batmal}, Theorem $\dis{1.1}$; \ci{jw}, Theorem $\dis{1.3}$)$\dis{}$\\
Inequality $\dis{\eqr{mttoda}}$ holds. In other words, the functional $\dis{J_\rho}$ is bounded from below on $\dis{H^1(\Si)^2}$ if and only if $\dis{\rho_i\le4\pi\min\lb1,1+\min_j\a_{ij}\rb}$ for $\dis{i=1,2}$.
\eth\

\bco\lab{mt}$\dis{}$\\
The functional $\dis{J_\rho}$ is coercive on $\dis{\ov H^1(\Si)^2}$ if and only if $\dis{\rho_i<4\pi\min\lb1,1+\min_j\a_{ij}\rb}$ for $\dis{i=1,2}$. If this occurs, then it has a global minimizer $\dis{u}$ which solves $\dis{\eqr{toda}}$.
\eco\

To overcome some difficulties of the problem, we will need a simple but essential topological result.\\

\ble$\dis{}$\\
Let $\dis{\Si}$ be a compact surface with $\dis{\chi(\Si)\le0}$. Then, there exist two curves $\dis{\g_1,\g_2}$, each of which is homeomorphic to a bouquet of $\dis{1+\ls\fr{-\chi(\Si)}2\rs}$ circles and two global projections $\dis{\Pi_i:\Si\to\g_i}$ such that (using the same notation as in Theorem $\dis{\ref{ex}}$):
\bit
\it $\dis{\g_1\cap\g_2=\es}$.
\it $\dis{p_{il}\in\g_i}$ for all $\dis{l\in\{1,\ds,L_i\}}$, $\dis{i=1,2}$.
\it $\dis{p_{0l}\nin\g_i}$ for all $\dis{l\in\{1,\ds,L_0\}}$, $\dis{i=1,2}$.
\eit
\ele\

\bfi[h]
\cen
\igr{./curve}
\cpt{The curves $\dis{\g_i}$}
\lab{fig1}
\efi

The proof of this lemma is quite intuitive and can be easily seen in Figure $\dis{\ref{fig1}}$.\\
If $\dis{\Si=\T^g}$ is a $\dis{g}$-torus, two retractions on disjoint bouquets of $\dis{g}$ circles can be easily built as for instance in \ci{bdm}. One can argue similarly with a connected sum $\dis{\Si=\P^{2k}}$ of an even number of copies of the projective plane, since this is homeomorphic to a connected sum of a $\dis{\T^{k-1}}$ and a Klein bottle, which in turn retracts on a circle; therefore, $\dis{\P^k}$ retracts on two disjoint bouquets of $\dis{k}$ circles. If instead $\dis{\Si}$ is a connected sum of an odd number of projective planes, one can argue as before setting the retractions constant on the last copy of $\dis{\P}$.\\
Notice that in all this case one has
$$g=1+\fr{-\chi(\T^g)}2,\q\q\q k=1+\fr{-\chi\lr\P^{2k}\rr}2=1+\ls\fr{-\chi\lr\P^{2k+1}\rr}2\rs$$
Finally, with a small deformation, the curves $\dis{\g_i}$ can be assumed to contain all the points $\dis{p_{il}}$ and they will not contain any of the other singular points. We can apply those deformations to $\dis{\g_1}$ without intersecting $\dis{\g_2}$ (or vice versa) because $\dis{\Si\bs\g_2}$ is pathwise connected.\\
The (non-weighted) barycenters on objects like $\dis{\g_i}$ have been considered in \ci{bdm} and their homology groups have been computed.\\

\bpr(\ci{bdm}, Proposition $\dis{3.2}$)\lab{homb}$\dis{}$\\
Let $\dis{\g}$ be a bouquet of $\dis{g}$ circles. Then, its barycenter sets verify
$$\wt H_q\lr\g^K\rr=\lb\bll\Z^{\bin{K+g-1}{g-1}}&\tx{if }q=2K-1\\0&\tx{if }q\ne2K-1\ea\ry.$$
\epr\

As mentioned in the introduction, we will have to deal with the set of weighted barycenters $\dis{\Si_{\rho,\a}}$ on $\dis{\Si}$ defined by $\dis{\eqr{wbar}}$. This is a subset of the space of the Radon measures $\dis{\mc M(\Si)}$ on $\dis{\Si}$ and it will be endowed with the $\dis{Lip'}$ norm, that is the norm of the dual space of Lipschitz functions:
$$\|\mu\|_{Lip'(\Si)}:=\sup_{\phi\in Lip(\Si),\|\phi\|_{Lip(\Si)}\le1}\lm\int_\Si\phi d\mu\rm.$$
We will denote as $\dis{d_{Lip'}}$ the corresponding distance.\\
The choice of this topology is somehow natural, since for any $\dis{x,y\in\Si}$ it holds $\dis{d_{Lip'}(\d_x,\d_y)\sim d(x,y)}$. Therefore, a copy of $\dis{\Si}$ is homeomorphically embedded in any $\dis{\Si_{\rho,\a}}$ and $\dis{\Si_1}$ is homeomorphic to $\dis{\Si}$.\\
When a measure is $\dis{Lip'}$-close to an element of $\dis{\Si_{\rho,\a}}$, it can be mapped onto this set, as proved in \cite{carmal}:\\

\ble(\ci{carmal}, Lemma $\dis{3.12}$)\lab{ret}$\dis{}$\\
For any $\dis{\rho\in\R}$, $\dis{\a=(\a_1,\ds,\a_J)}$ there exist $\dis{\e_0>0}$ and a continuous retraction
$$\psi_{\rho,\a}:\{\mu\in\mc M(\Si);\,d_{Lip'}(\mu,\Si_{\rho,\a})<\e_0\}\to\Si_{\rho,\a}$$
In particular, if $\dis{\mu_n\us{n\to\ity}\wk\s}$ for some $\dis{\s\in\Si_{\rho,\a}}$, then $\dis{\psi_{\rho,\a}(\mu_n)\us{n\to\ity}\to\s}$.
\ele\

At some point we will be under the assumptions of Lemma $\dis{\ref{ret}}$ for the function
\beq\lab{fui}
f_{i,u}:=\fr{\wt h_ie^{u_i}}{\int_\Si\wt h_ie^{u_i}dV_g},
\eeq
for one or both $\dis{i=1,2}$.\\
To express this alternative it will be useful to introduce the topological join of two sets $\dis{X}$ and $\dis{Y}$, that is
$$X\st Y:=\fr{X\x Y\x[0,1]}\sim$$
where $\dis{\sim}$ is the identification given by
$$(x,y,0)\sim(x,y',0)\q\fa\,x\in X,\,\fa\,y,y'\in Y,\q\q\q(x,y,1)\sim(x',y,1)\q\fa\,x,x'\in X,\,\fa\,y\in Y$$
The elements of $\dis{X\st Y}$ will be denoted by the formal expression $\dis{(1-t)x+ty}$.\\
The homology groups of the topological join can be easily expressed by the homology groups of $\dis{X}$ and $\dis{Y}$, since the $\dis{X\st Y}$ is the smash product of $\dis{X}$, $\dis{Y}$ and a copy of $\dis{\S^1}$. For details about smash products see for instance \ci{hat}.\\

\bpr (\ci{hat}, Theorem $\dis{3.21})$\lab{join}$\dis{}$\\
It holds
$$\wt H_q(X\st Y)=\Pl_{q'=0}^q\wt H_{q'}(X)\ox\wt H_{q-q'-1}(Y)$$
In particular, one has
$$\wt\b_q(X\st Y)=\sum_{q'=0}^q\wt\b_{q'}(X)\wt\b_{q-q'-1}(Y)$$
and
$$\sum_{q=0}^{+\ity}\wt\b_q(X\st Y)=\sum_{q'=0}^{+\ity}\wt\b_{q'}(X)\sum_{q''=0}^{+\ity}\wt\b_{q''}(Y).$$
\epr\

In particular, we will consider the topological join
$$\g_{\st,\rho,\wt\a}:=(\g_1)_{\rho_1,\wt\a_1}\st(\g_2)_{\rho_2,\wt\a_2}$$
between the two weighted barycenter sets defined on the curves $\dis{\g_1}$ with the multi-indices $\dis{\wt\a_i:=\lr\wt\a_{i1},\ds,\wt\a_{iL_i}\rr}$.\\

Let us now report the compactness result for $\dis{\eqr{toda}}$ from \ci{lwz}. We first introduce a finite set of couple of numbers, which represent the possible local blow-up values in a singular point $\dis{p}$ with coefficients $\dis{\a_1=\a_1(p)}$, $\dis{\a_2=\a_2(p)}$.\\

\bde\lab{lambda}$\dis{}$\\
For a couple of numbers $\dis{(\a_1,\a_2)}$ which are both greater than $\dis{-1}$, we set $\dis{\G_{\a_1,\a_2}\subset\mb R^2}$ as the piece of ellipse defined by
$$\G_{\a_1,\a_2}=\lb(\s_1,\s_2):\s_1^2-\s_1\s_2+\s_2^2-4\pi(1+\a_1)\s_2-4\pi(1+\a_2)\s_2=0\rb$$
We then define iteratively $\dis{\L_{\a_1,\a_2}\sub\G_{\a_1,\a_2}}$ via the following rules:
\bit
\it $\dis{\L_{\a_1,\a_2}}$ contains the points
$$(4\pi(1+\a_1),0),\q(0,4\pi(1+\a_2)),\q(4\pi(1+\a_1),4\pi(2+\a_1+\a_2))$$
$$(4\pi(2+\a_1+\a_2),4\pi(1+\a_2)),\q(4\pi(2+\a_1+\a_2),4\pi(2+\a_1+\a_2)).$$
\it If $\dis{(\s_1,\s_2)\in\L_{\a_1,\a_2}}$, then any $\dis{(\s'_1,\s'_2)\in\G_{\a_1,\a_2}}$ with $\dis{\s'_1=\s_1+4\pi n}$ for some $\dis{n\in\N}$ and $\dis{\s'_2\ge\s_2}$ belongs to $\dis{\L_{\a_1,\a_2}}$.
\it If $\dis{(\s_1,\s_2)\in\L_{\a_1,\a_2}}$, then any $\dis{(\s'_1,\s'_2)\in\G_{\a_1,\a_2}}$ with $\dis{\s'_2=\s_2+4\pi n}$ for some $\dis{n\in\N}$ and $\dis{\s'_1\ge\s_1}$ belongs to $\dis{\L_{\a_1,\a_2}}$.
\eit
\ede\

\bde$\dis{}$\\
Given $\dis{\L_{\a_1,\a_2}}$ as in Definition $\dis{\ref{lambda}}$, we define
$$\L_i:=\lb 4\pi n+\sum_{j\in\mc I}\s_j;\,n\in\N,\,\mc I\sub\{1,\ds,J\},\,\s_j\in\pi_i(\lr\G_{\a_{1j,\a_{2j}}}\rr)\rb,\q i=1,2,$$
where $\dis{\pi_i}$ is the projection on the $\dis{i}$-th component of $\dis{\R^2}$, and we set
$$\L:=(\L_1\x\R)\cup(\R\x\L_2)$$
\ede\

From the blow-up quantization in \ci{lwz} and an argument from \ci{batman,bremer} one finds a global compactness result.\\

\bth\lab{comp}$\dis{}$\\
If $\dis{\rho}$ belongs to a fixed compact set of $\dis{\R_+^2\bs\L}$, then the family of solutions of $\dis{\eqr{toda}}$ on $\dis{\ov H^1(\Si)}$ is uniformly bounded in $\dis{W^{2,p}(\Si)}$ for some $\dis{p>1}$.
\eth\

When dealing with compactness, there will be an essential result from \ci{luc} which helps to bypass the issue of Palais-Smale condition (which is not known to hold either for $\dis{\eqr{toda}}$ or for $\dis{\eqr{liou}}$). It basically states the Palais-Smale condition holds for some sequences on a dense set of $\dis{\rho}$ (see also \ci{djlw,str}); combining with Theorem $\dis{\ref{comp}}$ we deduce:\\

\ble\lab{def}$\dis{}$\\
Suppose $\dis{\rho\nin\L}$ and let $\dis{a<b\in\R}$ be such that $\dis{\eqr{toda}}$ has no solutions satisfying $\dis{a\le J_\rho\le b}$. Then, $\dis{J_\rho^a}$ is a deformation retract of $\dis{J_\rho^b}$.
\ele\

Moreover, compactness of solutions implies that $\dis{J_\rho}$ is bounded from above in its critical points, therefore we have the following:\\

\bco\lab{contr}$\dis{}$\\
Suppose $\dis{\rho\nin\L}$. Then, there exists $\dis{L>0}$ such that $\dis{J_\rho^L}$ is a deformation retract of $\dis{H^1(\Si)^2}$. In particular, it is contractible.
\eco\

Lemma $\dis{\ref{def}}$ and Corollary $\dis{\ref{contr}}$ allow to write in a simpler form the usual (weak) Morse inequalities $\dis{\mc C_q(a,b)\ge\b_q\lr J_\rho^b,J_\rho^a\rr}$.\\

\ble\lab{num}$\dis{}$\\
Suppose $\dis{\rho\nin\L}$ and $\dis{J_\rho}$ is a Morse function. Then, there exists $\dis{L>0}$ such that $\dis{\mc C_q(-L,L)\ge\wt\b_q\lr J_\rho^{-L}\rr.}$
In particular,
$$\#\tx{Solutions of }\eqr{toda}=\sum_{q=0}^{+\ity}\mc C_q\ge\sum_{q=0}^{+\ity}\mc C_q(-L,L)\ge\sum_{q=0}^{+\ity}\wt\b_q\lr J_\rho^{-L}\rr.$$
\ele\

Finally, the density results in \ci{bdm,dem} (which in turn use a transversality theorem from \ci{st}) can be immediately adapted from the scalar case to our purposes.\\

\bth\lab{dense}$\dis{}$\\
There exists an open dense set of $\dis{D\subset\mc M^2(\Si)\x L^\ity(\Si)^2}$ such that for any ${(g,h_1,h_2)\in D}$ $\dis{J_\rho}$ is a Morse function.
\eth\

\sec{Mapping weighted barycenters into sub-levels of $\dis{J_\rho}$}\

The aim of this section is to build a map which sends elements of $\dis{\g_{\st,\rho,\wt\a}}$ into arbitrarily low sub-levels of $\dis{J_\rho}$. We will actually build a family of maps $\dis{\Phi_\l}$ depending on a positive parameter $\dis{\l}$ such that $\dis{J_\rho\c\Phi_\l}$ attains negative values which are arbitrarily large in absolute value as $\dis{\l}$ gets larger.\\
The map $\dis{\Phi_\l}$ will be built starting from the standard bubbles
$$\ph_{\a,\l,p}=\log\lr\fr{\l^{1+\a}}{1+(\l d(\cd,p))^{2(1+\a)}}\rr^2,$$
which arise in study of $\dis{\eqr{liou}}$ under several aspects such as blow-up, compactness and characterization of sub-levels; functions with similar properties have been introduced for the case of $\dis{\eqr{toda}}$ in \ci{bjmr,batmal,jw,mn}.\\
There are basically two difficulties in building bubble-like functions which depend continuously on elements in $\dis{\g_{\st,\rho,\wt\a}}$: first of all, in the presence of a singular point $\dis{p_{il}}$ the parameter $\dis{\a}$ cannot switch suddenly from $\dis{0}$ to $\dis{\wt\a_{il}}$; moreover, we must be very careful of what happens in the two endpoints of the join, that is when one of the $\dis{(\g_i)_{\rho_i,\wt\a_i}}$ is identified to one point. However, we are able to fix both problems by arguing as in \ci{carmal} to avoid issues with singular points and as in \ci{bjmr} for what concerns the endpoints of the join.\\

\bth\lab{test}$\dis{}$\\
Given
$$\s_i=\sum_{x_{ik}\in\mc J_i}t_{ik}\d_{x_{ik}}\in(\g_i)_{\rho_i,\wt\a_i}\q\q\q\tx{for }i=1,2,\q\q\q\z=(1-t)\s_1+t\s_2\in\g_{\st,\rho,\wt\a}$$
define, for $\dis{\l>0}$,
$$\Phi_\l(\z)=\ph_{\l,\z}=(\ph_{1,\l,\z},\ph_{2,\l,\z})=\lr v_1-\fr{v_2}2,v_2-\fr{v_1}2\rr,$$
where
$$\l_{1,t}=\l(1-t),\q\q\q\l_{2,t}=\l t$$
\beq\lab{delta}
\d=\min\lb\min_{i=1,2,l=1,\ds,L_0}d(\g_i,p_{0l}),\fr{\min_{i=1,2,l\ne l'=1,\ds,L_i}d(p_{il},p_{il'})}2\rb
\eeq
\beq\lab{beta}
\b_{ik}=\lb\bll0&\tx{if }d:=\min_ld\lr x_{ik},p_{il}\rr\ge\d\\\fr{\wt\a_{i\wt l}\log\fr{\d}d}{\log\max\{2,\l_{i,t}\}-\wt\a_{i\wt l}\log\fr{\d}d}&\tx {if }d=d\lr x_{ik},p_{i\wt l}\rr\in\ls\max\{2,\l_{i,t}\}^{-\fr{1}{1+\wt\a_{i\wt l}}}\d,\d\rr\\\wt\a_{il}&\tx{if }d<\max\{2,\l_{i,t}\}^{-\fr{1}{1+\wt\a_{i\wt l}}}\d\ea\ry.
\eeq
and
\beq\lab{vi}
v_i=v_{i,\l,\z}=\log\sum_{x_{ik}\in\mc J_i}\fr{t_{ik}}{\lr1+\l_{i,t}^2d(\cd,x_{ik})^{2(1+\b_{ik})}\rr^2}.
\eeq
Then,
$$J_\rho(\ph_{\l,\z})\us{\l\to+\ity}\to-\ity\q\q\tx{uniformly for }\z\in\g_{\st,\rho,\wt\a}$$
\eth\

Notice that $\dis{\d}$ has been taken small enough so that the $\dis{p_{i\wt l}}$ which minimizes the distance between the singular points is uniquely determined if this distance is less than $\dis{\d}$. Furthermore, $\dis{\Phi_\l}$ is well defined because when $\dis{t=0}$, $\dis{v_1}$ vanishes and $\dis{v_2}$ depends only on $\dis{(\g_2)_{\rho_2,\a_2}}$, and the same occurs when $\dis{t=1}$ exchanging the roles of $\dis{v_1}$ and $\dis{v_2}$.\\
The choice of $\dis{\b_{ik}}$ has actually been made in such a way that, in the neck regions, it verifies
\beq\lab{deltad}
\max\{2,\l_{i,t}\}^\fr{\b_{ik}}{\a_{i\wt l}\lr1+\b_{ik}\rr}=\fr{\d}{d}.
\eeq
The proof of Theorem $\dis{\ref{test}}$ will follow by giving separate estimates for the three parts of $\dis{J_\rho}$, which will be provided respectively in the three following lemmas.\\

\ble\lab{q}$\dis{}$\\
Let $\dis{\z}$, $\dis{\ph_{\l,\z}}$ be as in Theorem $\dis{\ref{test}}$. Then,
$$\int_\Si Q(\ph_{\l,\z})\le8\pi\chi_{\wt\a_1}(\mc J_1)\log\max\{1,\l_{1,t}\}+8\pi\chi_{\wt\a_2}(\mc J_2)\log\max\{1,\l_{2,t}\}+C.$$
\ele\

\ble\lab{mean}$\dis{}$\\
Let $\dis{\z}$, $\dis{\ph_{\l,\z}}$ be as in Theorem $\dis{\ref{test}}$. Then,
$$\ov{\ph_{i,\l,\z}}=-4\log\max\{1,\l_{i,t}\}+2\log\max\{1,\l_{3-i,t}\}+O(1).$$
\ele\

\ble\lab{log}$\dis{}$\\
Let $\dis{\z}$, $\dis{\ph_{\l,\z}}$ be as in Theorem $\dis{\ref{test}}$. Then,
$$\log\int_\Si\wt h_ie^{\ph_{i,\l,\z}}dV_g=-2\log\max\{1,\l_{i,t}\}+2\log\max\{1,\l_{3-i,t}\}+O(1).$$
\ele\

Notice that, to prove Theorem $\dis{\ref{test}}$ we only need an estimate from below in Lemmas $\dis{\ref{mean}}$ and $\dis{\ref{log}}$. However, in Lemma $\dis{\ref{mean}}$ the same argument for the proof gives also an upper bound, whereas the estimate from above in Lemma $\dis{\ref{log}}$ will be needed later in this work.\\
We will only show the proof of Lemma $\dis{\ref{q}}$, since the others follow closely the proof given in \ci{car}, Proposition $\dis{4.1}$ and \ci{bjmr}, Proposition $\dis{3.3}$.\\

\bpf[Proof of Lemma $\dis{\ref{q}}$]$\dis{}$\\
First of all, we notice that
\beq\lab{qphi}
Q(\ph_{\l,\z})=\fr{1}4\lr|\n v_1|^2-\n v_1\cd\n v_2+|\n v_2|^2\rr.
\eeq
Since it holds
$$\n v_i=\fr{\sum_k\fr{-2(1+\b_{ik})t_{ik}\l_{i,t}^2d(\cd,x_{ik})^{1+2\b_{ik}}\n d(\cd,x_{ik})}{\lr1+\l_{i,t}^2d(\cd,x_{ik})^{2\lr1+\b_{ik}\rr}\rr^3}}{\sum_k\fr{t_{ik}}{\lr1+\l_{i,t}^2d(\cd,x_{ik})^{2\lr1+\b_{ik}\rr}\rr^2}}$$
and $\dis{|\n d(\cd,x_{ik})|=1}$ almost everywhere, we find
\bea
\nn|\n v_i|&\le&\fr{\sum_k\fr{4(1+\b_{ik})t_{ik}\l_{i,t}^{2\lr1+\b_{ik}\rr}d(\cd,x_{ik})^{1+2\b_{ik}}}{\lr1+\l_{i,t}^2d(\cd,x_{ik})^{2\lr1+\b_{ik}\rr}\rr^3}}{\sum_k\fr{t_{ik}}{\lr1+\l_{i,t}^2d\lr(\cd,x_{ik})\rr^{2\lr1+\b_{ik}\rr}\rr^2}}\\
\lab{nvi}&\le&\max_k\ub{\fr{4(1+\b_{ik})\l_{i,t}^2d(\cd,x_{ik})^{1+2\b_{ik}}}{1+\l_{i,t}^2(d(\cd,x_{ik}))^{2\lr1+\b_{ik}\rr}}}_{=:m_{ik}}.
\eea
In view of these estimates, we divide $\dis{\Si}$ into a finite number of regions depending on which of the $\dis{m_{1k}}$'s attains the maximum:
$$A_{1k}:=\lb x\in\Si:m_{1k}(x)=\max_{k'}m_{1k'}(x)\rb.$$
Similarly, we will define, in dependence of the $\dis{m_{2k}}$'s:
$$A_{2k}:=\lb x\in\Si:m_{2k}(x)=\max_{k'}m_{2k'}(x)\rb.$$
Moreover, we can easily see that the following estimates hold for $\dis{m_{ik}}$:
\beq\lab{mk}
m_{ik}\le\lb\bl\fr{4\lr1+\b_{ik}\rr}{d\lr\cd,x_{ik}\rr}\\4(1+\b_{ik})\l_{i,t}^2d(\cd,x_{ik})^{1+2\b_{ik}}\ea\ry.
\eeq
We will estimate the mixed term first. Basically, since the points $\dis{x_{ik}}$ belong to $\dis{\g_i}$ and the curves $\dis{\g_i}$'s are disjoint, we only have summable singularities and therefore the integral of $\dis{\n v_1\cd\n v_2}$ is uniformly bounded.\\
Therefore, from $\dis{\eqr{nvi}}$ and the first inequality in $\dis{\eqr{mk}}$, one finds
\bey
\lm\int_\Si\n v_1\cd\n v_2dV_g\rm&\le&\sum_{k,k'}\int_{A_{1k}\cap A_{2k'}}|\n v_1||\n v_2|dV_g\le\\
&&\sum_{k,k'}16(1+\b_{1k})(1+\b_{2k'})\int_{A_{1k}\cap A_{2k'}}\fr{dV_g}{d(\cd,x_{1k})d(\cd,x_{2k'})}
\eey
We then notice that, by the definition of $\dis{\eqr{delta}}$, the distance between $\dis{\g_1}$ and $\dis{\g_2}$ is at least $\dis{2\d}$, so $\dis{B_\d(x_{1k})\cap B_\d(x_{2k'})=\es}$ for any choice of $\dis{k,k'}$. Therefore,
\bey
\int_{A_{1k}\cap A_{2k'}}\fr{dV_g}{d(\cd,x_{1k})d(\cd,x_{2k'})}&\le&\int_{A_{1k}\cap A_{2k'}\bs B_\d\lr x_{1k}\rr}\fr{dV_g}{\d d(\cd,x_{2k'})}\\
&+&\int_{A_{1k}\cap A_{2k'}\bs B_\d\lr x_{2k'}\rr}\fr{dV_g}{\d d(\cd,x_{1k})}\\
&\le&\fr{1}\d\int_\Si\lr\fr{1}{d(\cd,x_{2k'})}+\fr{1}{d(\cd,x_{1k})}\rr dV_g\\
&\le&C_\d,
\eey
hence, being the number of $\dis{k,k'}$ bounded from above depending on $\dis{\rho}$ and $\dis{\wt\a_{il}}$'s only, we obtain
\beq\lab{v1v2}
\lm\int_\Si\n v_1\cd\n v_2dV_g\rm\le C.
\eeq
Now, we consider the term involving $\dis{|\n v_1|^2}$. We split the integral into the sets $\dis{A_{1k}}$ defined above.
\bea
\nn\int_\Si|\n v_1|^2dV_g&\le&\sum_k\int_{A_{1k}}m_{1k}^2dV_g\\
\lab{v1sum}&\le&\sum_k\lr\int_{\Si\bs B_{\l_{1,t}^{-\fr{1}{1+\b_{1k}}}}(x_{1k})}m_{1k}^2dV_g\ry+\ly\int_{B_{\l_{1,t}^{-\fr{1}{1+\b_{1k}}}}(x_{1k})}m_{1k}^2dV_g\rr
\eea
Outside the balls we will apply the first estimate in $\dis{\eqr{nvi}}$:
\bea
\nn\int_{\Si\bs B_{\l_{1,t}^{-\fr{1}{1+\b_{1k}}}}(x_{1k})}|\n v_1|^2dV_g&\le&16(1+\b_{1k})^2\int_{\Si\bs B_{\l_{1,t}^{-\fr{1}{1+\b_{1k}}}}(x_{1k})}\fr{dV_g}{d(\cd,x_{1k})^2}\\
\nn&\le&32\pi(1+\b_{1k})^2\log\max\lb1,\l_{1,t}^{\fr{1}{1+\b_{1k}}}\rb+C\\
\lab{v1out}&\le&32\pi(1+\b_{1k})\log\max\{1,\l_{1,t}\}+C.
\eea
The integral inside the balls is actually uniformly bounded, as can be seen using now the second estimate in $\dis{\eqr{nvi}}$:
\bea
\nn\int_{B_{\l_{1,t}^{-\fr{1}{1+\b_{1k}}}}\lr x_{1k}\rr}|\n v_1|^2dV_g&\le&16(1+\b_{1k})^2\l_{1,t}^4\int_{B_{\l_{1,t}^{-\fr{1}{1+\b_{1k}}}}\lr x_{1k}\rr}d(\cd,x_{1k})^{2(1+2\b_{1k})}dV_g\\
\nn&\le&C_{\b_{1k}}\l_{1,t}^4\lr\l_{1,t}^{-\fr{1}{1+\b_{1k}}}\rr^{4(1+\b_{1k})}\\
\lab{v1in}&\le&C
\eea
Observing that, from the definitions of $\dis{\eqr{weight}}$ and $\dis{\eqr{beta}}$, one has
$$\sum_k(1+\b_{1k})\le\chi_{\wt\a_1}(\mc J_1),$$
one can now deduce from $\dis{\eqr{v1sum}}$, $\dis{\eqr{v1out}}$, $\dis{\eqr{v1in}}$:
\beq\lab{v1}
\int_\Si|\n v_1|^2\le32\pi\chi_{\wt\a_1}(\mc J_1)\log\min\{C,\l_{1,t}\}+C.
\eeq
The same argument gives a similar estimate for $\dis{\int_\Si|\n v_2|^2}$, therefore putting together $\dis{\eqr{v1}}$ with $\dis{\eqr{v1v2}}$ and $\dis{\eqr{qphi}}$ we get the conclusion.
\epf\

\bpf[Proof of Theorem $\dis{\ref{test}}$]$\dis{}$\\
Since on $\dis{(\g_i)_{\rho_i,\wt\a_i}}$ one has $\dis{\rho_i<4\pi\chi_{\wt\a_i}(\mc J_i)}$, and moreover $\dis{\max\{\l_{1,t},\l_{2,t}\}\ge\fr{\l}2}$, Lemmas $\dis{\ref{q}}$, $\dis{\ref{mean}}$ and $\dis{\ref{log}}$ yield
\bey
J_\rho(\ph_{\l,\z})&=&\int_\Si Q(\ph_{\l,\z})dV_g+\sum_{i=1}^2\rho_i\lr\ov{\ph_{i,\l,\z}}-\log\int_\Si\wt h_ie^{\ph_{i,\l,\z}}dV_g\rr\\
&\le&\lr8\pi\chi_{\wt\a_1}(\mc J_1)-2\rho_1\rr\log\max\{1,\l_{1,t}\}+(8\pi\chi_{\wt\a_2}(\mc J_2)-\rho_2)\log\max\{1,\l_{2,t}\}+C\\
&\le&\max\{8\pi\chi_{\wt\a_1}(\mc J_1)-2\rho_1,8\pi\chi_{\wt\a_2}(\mc J_2)-2\rho_2\}\log\max\{1,\l_{1,t},\l_{2,t}\}+C\\
&\le&\max\{8\pi\chi_{\wt\a_1}(\mc J_1)-2\rho_1,8\pi\chi_{\wt\a_2}(\mc J_2)-2\rho_2\}\log\max\{1,\l\}+C\\
&\us{\l\to+\ity}\to&-\ity
\eey
uniformly in $\dis{\z\in\g_{\st,\rho,\wt\a}}$, which is what we wished to prove.
\epf\

\sec{Analysis of sub-levels and improved Moser-Trudinger inequalities}\

In this section we are going to provide information on the sub-levels $\dis{J_\rho^{-L}}$.\\
A key point is the so-called improved Moser-Trudinger inequality. Basically, we show that under certain conditions of the spreading on $\dis{u_1}$ and $\dis{u_2}$ the constant in Moser-Trudinger inequality can be improved and from this fact we deduce information about the low sub-levels of $\dis{J_\rho}$.\\
The idea was introduced by Chen and Li \ci{cl01} for the Liouville equation and extended in \ci{dja,mal08}, in \ci{carmal} for the singular case and in \ci{bjmr,mn,mr13} for the Toda system. Some results presented in this section will be adapted from the aforementioned papers, so their proof will be skipped.\\
The main result of this section is the following:\\

\bth\lab{deps}$\dis{}$\\
Suppose $\dis{\rho\nin\L}$. Then, for any $\dis{\e>0}$ there exists $\dis{L=L(\e)>0}$ such that any $\dis{u\in J_\rho^{-L}}$ verifies, for some $\dis{i=1,2}$,
$$d_{Lip'}(f_{i,u},\Si_{\rho_i,\wt\a_i})<\e,$$
where $\dis{f_{i,u}}$ is defined by $\dis{\eqr{fui}}$.
\eth\

To adapt the original argument to the case of Toda system we first need a covering lemma (\ci{bjmr}, Lemma $\dis{4.1}$; \ci{mn}, Lemma $\dis{3.2}$; \ci{mr13}, Lemma $\dis{2.5}$).\\
With respect to the previous works, we have to take into account the singularities and consider sets which contain at most one negative singularity. Anyway, the proof can be adapted step-by-step with the conditions on the singular points still holding, so we will omit it.\\

\ble\lab{cov}$\dis{}$\\
Let $\dis{\d>0,\t>0,M_1,N_1,M_2,N_2\in\N}$ be given numbers, $\dis{\{l_1,\ds,l_{M_i}\}\sub\{1,\ds,L_i\}}$ selections of indices, $\dis{f_1,f_2\in L^1(\Si)}$ be non-negative functions with $\dis{\int_\Si f_idV_g=1}$ and $\dis{\{\O_{im}\}_{i=1,2}^{m=1,\ds,M_i+N_i}}$ be measurable subsets of $\dis{\Si}$ such that
$$d(\O_{im},\O_{im'})\ge\d\q\q\q\fa\,i=1,2,\,\fa\,m,m'=1,\ds,M_i+N_i,\,m\ne m'$$
$$d(p_{il},\O_{im})\ge\d\q\q\q\fa\,i=1,2,\,\fa\,m=1,\ds,M_i+N_i,\,\fa\,l=1,\ds,L_i,\,l\ne l_m$$
$$\int_{\O_{im}}f_idV_g\ge\t\q\q\q\fa\,i=1,2,\,\fa\,m=1,\ds,M_i+N_i.$$
Then, there exist $\dis{\wt\d>0,\wt\t>0}$, independent of $\dis{f_1,f_2}$, and $\dis{\{\O_m\}_{m=1,\ds,\max_i\{M_i+N_i\}}}$ such that
$$d(\O_m,\O_{m'})\ge\wt\d\q\q\q\fa\,m,m'=1,\ds,\max_{i=1,2}\{M_i+N_i\},\,m\ne m'$$
$$d(p_{il},\O_m)\ge\wt\d\q\q\q\fa\,i=1,2,\,\fa\,m=1,\ds,\max_{i=1,2}\{M_i+N_i\},\,\fa\,l=1,\ds,L_i,\,l\ne l_m$$
$$\int_{\O_m}f_idV_g\ge\wt\t\q\q\q\fa\,i=1,2,\,\fa\,m=1,\ds,\max_{i=1,2}\{M_i+N_i\}.$$
\ele\

The next lemma is what is usually called an improved Moser-Trudinger inequality.\\
It essentially states that if both $\dis{u_1}$ and $\dis{u_2}$ are spread in sets which contain at most one singular point, then the constant $\dis{4\pi}$ in $\dis{\eqr{mttoda}}$ can be multiplied by a number depending on how many these sets are and on the singular points they contain.\\

\ble\lab{impr}$\dis{}$\\
Let $\dis{\d>0,\t>0,M_1,N_1,M_2,N_2\in\N}$ be given numbers, $\dis{\{l_1,\ds,l_{M_i}\}\sub\{1,\ds,L_i\}}$ selections of indices and $\dis{\{\O_{im}\}_{i=1,2}^{m=1,\ds,M_i+N_i}}$ be measurable subsets of $\dis{\Si}$ such that
$$d(\O_{im},\O_{im'})\ge\d\q\q\q\fa\,i=1,2,\,\fa\,m,m'=1,\ds,M_i+N_i,\,m\ne m'$$
$$d(p_{il},\O_{im})\ge\d\q\q\q\fa\,i=1,2,\,\fa\,m=1,\ds,M_i+N_i,\,\fa\,l=1,\ds,L_i,\,l\ne l_m$$
Then, for any $\dis{\e>0}$ there exists $\dis{C>0}$, not depending on $\dis{u}$, such that any $\dis{u=(u_1,u_2)\in H^1(\Si)^2}$ satisfying
$$\int_{\O_{im}}f_{i,u}dV_g\ge\t\q\q\q\fa\,i=1,2,\,\fa\,m=1,\ds,M_i+N_i$$
verifies
$$\sum_{i=1}^2\lr N_i+\sum_{m=1}^{M_i}(1+\wt\a_{il_m})\rr\log\int_\Si\wt h_ie^{u_i-\ov{u_i}}dV_g\le\fr{1+\e}{4\pi}\int_\Si Q(u)dV_g+C.$$
\ele\

\bpf$\dis{}$\\
It will not be restrictive to suppose $\dis{M_1+N_1\ge M_2+N_2}$.\\
We apply Lemma $\dis{\ref{cov}}$ with $\dis{f_i=f_{i,u}}$ and we get a family of sets $\dis{\{\O_m\}_{m=1}^{M_1+N_1}}$ satisfying
$$d(\O_m,\O_{m'})\ge\wt\d>0\q\fa\,m\ne m',\q\q\q\int_{\O_m}f_{i,u}\ge\wt\t>0\q\fa\,m=1,\ds,M_i+N_i$$
Let us now consider, for any $\dis{m=1,\ds,M_1+N_1}$, the cut-off function $\dis{g_m:=\max\lb0,1-\fr{2}{\wt\d}d(\cd,\O_m)\rb}$; it verifies
$$\mf 1_{\O_m}\le g_m\le\mf 1_{\wt\O_m},\q\q\q|\n g_m|\le C_{\wt\d,\Si}\mf1_{\wt\O_m}\q\q\q\tx{with }\wt\O_m=B_{\fr{\wt\d}2}(\O_m).$$
We now take $\dis{v_i\in L^\ity(\Si)}$ with $\dis{\ov{v_i}=0}$ and we set $\dis{w_i:=u_i-v_i-\ov{u_i}}$ (which will also have null average). Therefore, we find
\bea
\nn\log\int_\Si\wt h_1e^{u_i-\ov{u_i}}dV_g&\le&\log\lr\fr{1}{\wt\t}\int_{\O_m}\wt h_ie^{u_i-\ov{u_i}}dV_g\rr\\
\nn&=&\log\lr\fr{1}{\wt\t}\int_{\O_m}\wt h_ie^{v_i+w_i}dV_g\rr\\
\nn&\le&\log\fr{1}{\wt\t}+\|v_i\|_{L^\ity(\O_m)}+\log\int_{\O_m}\wt h_ie^{w_i}dV_g\\
\lab{imprlog}&\le&\log\fr{1}{\wt\t}+\|v_i\|_{L^\ity(\Si)}+\log\int_\Si\wt h_ie^{g_mw_i}dV_g.
\eea
Since $\dis{g_m\in Lip(\Si)}$, then $\dis{g_mw_i\in H^1(\Si)}$, so we can apply a Moser-Trudinger inequality on it.\\
To this purpose, we notice that, for any $\dis{\e>0}$,
\bey
\int_\Si|\n(g_mw_1)|^2dV_g&=&\int_\Si|g_m\n w_1+w_1\n g_m|^2dV_g\\
&=&\int_\Si\lr g_m^2|\n w_1|^2+2(g_m\n w_1)\cd(w_1\n g_m)+w_1^2|\n g_m|^2\rr dV_g\\
&\le&\int_\Si\lr(1+\e)g_m^2|\n w_1|^2+\lr1+\fr{1}\e\rr w_1^2|\n g_m|^2\rr dV_g\\
&\le&(1+\e)\int_{\wt\O_m}|\n w_1|^2dV_g+C_{\e,\wt\d,\Si}\int_{\wt\O_m}w_1^2dV_g.
\eey
In the same way, writing 
\beq\lab{xy}
\fr{1}3\lr|x|^2+x\cd y+|y|^2\rr=\fr{1}4|x|^2+\fr{1}{12}|x-2y|^2,
\eeq we get
\beq\lab{imprq}
\int_\Si Q(g_mw)dV_g\le(1+\e)\int_{\wt\O_m}Q(w)dV_g+C_{\e,\wt\d,\Si}\int_{\wt\O_m}\fr{1}3\lr w_1^2+w_1w_2+w_2^2\rr dV_g.
\eeq
At this point, we choose properly $\dis{w_i}$ (hence $\dis{v_i}$) in such a way to have a control of its $\dis{L^2}$ norm. Taking an orthonormal frame $\dis{\{\ph_n\}_{n=1}^\ity}$ of eigenfunctions for $\dis{-\D}$ on $\dis{\ov H^1(\Si)}$ with a non-decreasing sequence of associated positive eigenvalues $\dis{\{\l_n\}_{n=1}^\ity}$ and writing $\dis{u_i=\ov{u_i}+\sum_{n=1}^\ity u_{in}\ph_n}$, we set $\dis{v_i=\sum_{n=1}^Nu_{in}\ph_n}$ for
$$N=N_{\e,\wt\d,\Si}:=\max\lb n\in\N:\l_n<\fr{C_{\e,\wt\d,\Si}}\e\rb.$$
This choice gives
$$
C_{\e,\wt\d,\Si}\int_\Si w_1^2dV_g\le\e\int_\Si|\n w_1|^2dV_g\le\e\int_\Si|\n u_1|^2dV_g
$$
and, through $\dis{\eqr{xy}}$,
\beq\lab{imprw}
C_{\e,\wt\d,\Si}\int_\Si\fr{1}3\lr w_1^2+w_1w_2+w_2^2\rr dV_g\le\e\int_\Si Q(w)dV_g\le\e\int_\Si Q(u)dV_g.
\eeq
Moreover, we get
\beq\lab{imprwi}
\ov{|w_i|}\le C_\Si\|\n w_i\|_{L^2(\Si)}\le \e\int_\Si Q(w)dV_g+C\le\e\int_\Si Q(u)dV_g+C
\eeq
and, since $\dis{v_i}$ belongs to a finite-dimensional space,
\beq\lab{imprvi}
\|v_i\|_{L^\ity(\Si)}\le C_N\|\n v_i\|_{L^2(\Si)}\le\e\int_\Si Q(v)dV_g+C\le\e\int_\Si Q(u)dV_g+C.
\eeq
Now, if $\dis{m=1,\ds,M_2+N_2}$, we apply the Moser-Trudinger inequality $\dis{\eqr{mttoda}}$ to $\dis{g_mw}$. Since these functions are supported on $\dis{\wt\O_m}$, we can replace $\dis{\wt h_i}$ by a smooth interpolation which is constant outside a neighborhood of $\dis{\wt\O_m}$: we take $\dis{\eta_m}$ satisfying
$$\eta_m(x):=\lb\bll1&\tx{if }x\in\wt\O_m\\0&\tx{if }x\nin B_{\fr{\wt\d}4\lr\wt\O_m\rr}\ea\ry\q\q\q\wt h_{im}:=\eta_m\wt h_i+1-\eta_m=\lb\bll\wt h_i&\tx{if }x\in\wt\O_m\\1&\tx{if }x\nin B_{\fr{\wt\d}4\lr\wt\O_m\rr}\ea\ry.$$
In this way, we can consider only the singularities $\dis{p_{1l_m},p_{2l_m}}$ which lie inside $\dis{\O_m}$ (if there are any); from $\dis{\eqr{imprq}}$ and $\dis{\eqr{imprwi}}$ we get
\bey
\sum_{i=1}^2(1+\wt\a_{il_m})\log\int_\Si\wt h_ie^{g_mw_i}dV_g&=&\sum_{i=1}^2(1+\wt\a_{il_m})\log\int_\Si\wt h_{im}e^{g_mw_i}dV_g\\
&\le&\sum_{i=1}^2(1+\wt\a_{il_m})\ov{g_mw_i}+\fr{1}{4\pi}\int_\Si Q(g_mw)dV_g+C\\
&\le&\sum_{i=1}^2(1+\wt\a_{il_m})\lr\|g_m\|_{L^\ity(\Si)}\ov{|w_i|}\rr+\fr{1+\e}{4\pi}\int_{\wt\O_m}Q(w)dV_g\\
&+&\fr{C_{\e,\wt\d,\Si}}{4\pi}\int_{\wt\O_m}\fr{1}3\lr w_1^2+w_1w_2+w_2^2\rr dV_g+C\\
&\le&\ov{|w_1|}+\ov{|w_2|}+\fr{1+\e}{4\pi}\int_{\wt\O_m}Q(w)dV_g\\
&+&\fr{C_{\e,\wt\d,\Si}}{4\pi}\int_{\wt\O_m}\fr{1}3\lr w_1^2+w_1w_2+w_2^2\rr dV_g+C\\
&\le&2\e\int_\Si Q(u)dV_g+\fr{1+\e}{4\pi}\int_{\wt\O_m}Q(w)dV_g\\
&+&\fr{C_{\e,\wt\d,\Si}}{4\pi}\int_{\wt\O_m}\fr{1}3\lr w_1^2+w_1w_2+w_2^2\rr dV_g+C.
\eey
Therefore, from $\dis{\eqr{imprlog}}$ and $\dis{\eqr{imprvi}}$ we deduce
\bea
\nn\sum_{i=1}^2(1+\wt\a_{il_m})\log\int_\Si\wt h_ie^{u_i-\ov{u_i}}dV_g&\le&\sum_{i=1}^2(1+\wt\a_{il_m})\lr\log\fr{1}{\wt\t}+\|v_i\|_{L^\ity(\Si)}+\log\int_\Si\wt h_ie^{g_mw_i}dV_g\rr\\
\nn&\le&\sum_{i=1}^2\lr\|v_i\|_{L^\ity(\Si)}+(1+\wt\a_{il_m})\log\int_\Si\wt h_ie^{g_mw_i}dV_g\rr+C\\
\nn&\le&3\e\int_\Si Q(u)dV_g+\fr{1+\e}{4\pi}\int_{\wt\O_m}Q(w)dV_g\\
\lab{imprsum}&+&\fr{C_{\e,\wt\d,\Si}}{4\pi}\int_{\wt\O_m}\fr{1}3\lr w_1^2+w_1w_2+w_2^2\rr dV_g+C.
\eea
For $\dis{m=M_2+N_2+1,\ds,M_1+N_1}$ we have estimates only for $\dis{u_1}$ on $\dis{\O_m}$, so we apply the scalar Moser-Trudinger inequality $\dis{\eqr{mtliou}}$. By $\dis{\eqr{xy}}$ we get the integral of $\dis{Q(g_mw)}$, then we argue as before.\\
Notice that if $\dis{m>M_i}$, then $\dis{p_{il_m}}$ is not defined so these calculations would not make sense, but in this case both the previous and the following calculations hold replacing $\dis{\wt\a_{il_m}}$ with $\dis{0}$.
\bey
(1+\wt\a_{1l_m})\log\int_\Si\wt h_1e^{g_mw_1}dV_g&=&(1+\wt\a_{1l_m})\log\int_\Si\wt h_{1m}e^{g_mw_1}dV_g\\
&\le&(1+\wt\a_{1l_m})\ov{g_mw_1}+\fr{1}{16\pi}\int_\Si|\n(g_mw_1)|^2dV_g+C\\
&\le&\ov{|w_1|}+\fr{1}{4\pi}\int_\Si Q(g_mw)dV_g+C\\
&\le&\e\int_\Si Q(u)dV_g+\fr{1+\e}{4\pi}\int_{\wt\O_m}Q(w)dV_g\\
&+&\fr{C_{\e,\wt\d,\Si}}{4\pi}\int_{\wt\O_m}\fr{1}3\lr w_1^2+w_1w_2+w_2^2\rr dV_g+C.
\eey
Then in this case we deduce
\bea
\nn(1+\wt\a_{1l_m})\log\int_\Si\wt h_1e^{u_1-\ov{u_1}}dV_g&\le&(1+\wt\a_{1l_m})\lr\log\fr{1}{\wt\t}+\|v_1\|_{L^\ity(\Si)}+\log\int_\Si\wt h_ie^{g_mw_1}dV_g\rr\\
\nn&\le&\|v_1\|_{L^\ity(\Si)}+(1+\wt\a_{1l_m})\log\int_\Si\wt h_1e^{g_mw_1}dV_g+C\\
\nn&\le&2\e\int_\Si Q(u)dV_g+\fr{1+\e}{4\pi}\int_{\wt\O_m}Q(w)dV_g\\
\lab{imprfin}&+&\fr{C_{\e,\wt\d,\Si}}{4\pi}\int_{\wt\O_m}\fr{1}3\lr w_1^2+w_1w_2+w_2^2\rr dV_g+C.
\eea
Finally, we sum together $\dis{\eqr{imprsum}}$ and $\dis{\eqr{imprfin}}$ for all the $\dis{m}$'s, exploiting $\dis{\eqr{imprw}}$ and the disjointness of the $\dis{\wt\O_m}$:
\bey
\sum_{i=1}^2\lr N_i+\sum_{m=1}^{M_i}(1+\wt\a_{il_m})\rr\log\int_\Si\wt h_ie^{u_i-\ov{u_i}}dV_g&=& \sum_{i=1}^2\sum_{m=1}^{M_2+N_2}(1+\wt\a_{il_m})\lr\log\int_\Si\wt h_ie^{u_i-\ov{u_i}}dV_g\rr\\
&+&\sum_{m=M_2+N_2+1}^{M_1+N_1}(1+\wt\a_{1l_m})\log\int_\Si\wt h_1e^{u_1-\ov{u_1}}dV_g\\
&\le&(2M_1+2M_1+M_2+N_2)\e\int_\Si Q(u)dV_g\\
&+&\fr{1+\e}{4\pi}\int_\Si Q(w)dV_g\\
&+&\fr{C_{\e,\wt\d,\Si}}{4\pi}\int_\Si\fr{1}3\lr w_1^2+w_1w_2+w_2^2\rr dV_g+C\\
&\le&(2M_1+2N_1+M_2+N_2)\e\int_\Si Q(u)dV_g\\
&+&\fr{1+2\e}{4\pi}\int_\Si Q(u)dV_g+C
\eey
which is, renaming $\dis{\e}$ properly, what we desired.
\epf\

Now we need another technical lemma, which relates the condition of spreading, needed for Lemma $\dis{\ref{impr}}$, and of concentration around a finite number of points.\\
Through this lemma, we can then use the improved Moser-Trudinger inequality to get information about the concentration which occurs on sub-levels $\dis{J_\rho^{-L}}$.\\
The following results will be extensions of the ones contained in \ci{bjmr,dm,mn,mr13} with suitable changes to take into account the singularities. Since the modifications are minimal, the proofs will be skipped.\\

\ble\lab{br}$\dis{}$\\
Let $\dis{i=1,2}$, $\dis{\chi_0>0}$, $\dis{\e,r>0}$ small enough, be such that
any $\dis{\mc J\sub\Si}$ satisfying $\dis{\chi_{\wt\a_i}(\mc J)\le\chi_0}$
verifies
$$\int_{\Cup_{x_k\in\mc J}B_r(x_k)}f_{i,u}dV_g<1-\e\q\q\q.$$
Then, there exist $\dis{\wt\e,\wt r>0}$, not depending on $\dis{u_i}$, $\dis{M,N\in\N}$, $\dis{\lb l_1,\ds,l_M\rb\sub\{1,\ds,L_i\}}$ and $\dis{\lb\wt x_m\rb_{m=1}^{M+N}}$ satisfying
$$N+\sum_{m=1}^M\lr1+\wt\a_{il_m}\rr>\chi_0,\q\q\q d(\wt x_m,p_l)\ge2\wt r\q\fa\,l\ne l_m$$
$$B_{2\wt r}(\wt x_m)\cap B_{2\wt r}(\wt x_{m'})=\es\q\fa\,m\ne m',\q\q\q\int_{B_{\wt r}(\wt x_m)}f_{i,u}dV_g\ge\wt\e\q\fa\,m=1,\ds,M+N.$$
\ele\

\ble\lab{intbr}$\dis{}$\\
For any $\dis{\e,r>0}$ small enough, there exists $\dis{L>0}$ such that,  if $\dis{u\in J_\rho^{-L}}$, then for at least one $\dis{i=1,2}$ there exists $\dis{\mc J\sub\Si}$ satisfying $\dis{4\pi\chi_{\wt\a_i}(\mc J_i)\le\rho_i}$ and
$$\int_{\Cup_{x_k\in\mc J_i}B_r(x_k)}f_{i,u}dV_g\ge1-\e.$$
\ele\

Now we have all the tools to prove Theorem $\dis{\ref{deps}}$.\\

\bpf[Proof of Theorem $\dis{\ref{deps}}$]$\dis{}$\\
It clearly suffices to prove the statement for small $\dis{\e}$.\\
We apply Lemma $\dis{\ref{intbr}}$ with $\dis{\e'=\fr{\e}4,r'=\fr{\e}2}$. It is not restrictive to suppose that the thesis of the lemma holds for $\dis{i=1}$, since the case $\dis{i=2}$ can be treated in the same way. Therefore, we get $\dis{\mc J\sub\Si}$, and we define
$$\s_u:=\sum_{x_k\in\mc J}t_k\d_{x_k}$$
where
$$t_k=\int_{B_{r'}(x_k)\bs\Cup_{k'=1}^{k-1}B_{r'}(x_{k'})}f_{1,u}dV_g+\fr{1}{|\mc J|}\int_{\Si\bs\Cup_{x_{k'}\in\mc J}B_{r'}(x_{k'})}f_{1,u}dV_g.$$
Notice that $\dis{\s_u\in\Si_{\rho_1,\wt\a_1}}$ because, from Lemma $\dis{\ref{intbr}}$ we find $\dis{\chi_{\wt\a_1}(\mc J)\le\rho_1}$ and the last inequality is actually strict because we are supposing $\dis{\rho\nin\L}$.\\
To conclude the proof it would suffice to get
\beq\lab{depsphi}
\lm\int_\Si\lr f_{1,u}-\s_u\rr\phi dV_g\rm\le\e\|\phi\|_{Lip(\Si)}\q\q\q\fa\,\phi\in Lip(\Si).
\eeq
In fact, following the definition of $\dis{d_{Lip'}}$, this would imply
\beq\lab{depslip}
d_{Lip'}(f_{1,u},\Si_{\rho_1,\wt\a_1})\le d_{Lip'}(f_{1,u},\s_u)=\sup_{\phi\in Lip(\Si),\|\phi\|_{Lip(\Si)}\le1}\lm\int_\Si\lr f_{1,u}-\s_u\rr\phi dV_g\rm<\e.
\eeq
We will divide the integral in $\dis{\eqr{depsphi}}$ into two points, studying separately what happens inside and outside the union of the $\dis{r'}$-balls centered at the points $\dis{x_m}$'s.\\
Outside the balls, for any $\dis{\phi\in Lip(\Si)}$ we have
\bea
\nn\lm\int_{\Si\bs\Cup_{x_k\in\mc J}B_{r'}(x_k)}(f_{1,u}-\s_u)\phi dV_g\rm&=&\lm\int_{\Si\bs\Cup_{x_k\in\mc J}B_{r'}(x_k)}f_{1,u}\phi dV_g\rm\\
\nn&\le&\|\phi\|_{L^\ity(\Si)}\int_{\Si\bs\Cup_{x_k\in\mc J}B_{r'}(x_k)}f_{1,u}dV_g\\
\nn&<&\e'\|\phi\|_{Lip(\Si)}\\
\lab{depsout}&=&\fr{\e}4\|\phi\|_{Lip(\Si)}.
\eea
On the other hand, we also find
\bea
\nn\lm\int_{\Cup_{x_k\in\mc J}B_{r'}(x_k)}(f_{1,u}-\s_u)\phi dV_g\rm&=&\lm\int_{\Cup_{x_k\in\mc J}B_{r'}(x_k)}f_{1,u}\phi dV_g\ry\\
\nn&-&\ly\sum_{x_k\in\mc J}\lr\int_{B_{r'}(x_k)\bs\Cup_{k'=1}^{k}B_{r'}(x_{k'})}f_{1,u}dV_g+\ry\ry\\
\nn&+&\ly\ly\fr{1}{|\mc J|}\int_{\Si\bs\Cup_{x_{k'}\in\mc J}B_{r'}(x_{k'})}f_{1,u}dV_g\rr\phi(x_k)\rm\\
\nn&=&\lm\sum_{x_k\in\mc J}\lr\int_{B_{r'}(x_k)\bs\Cup_{k'=1}^kB_{r'}(x_{k'})}f_{1,u}(\phi-\phi(x_k))dV_g\ry\ry\\
\nn&-&\ly\ly\fr{1}{|\mc J|}\int_{\Si\bs\Cup_{x_{k'}\in\mc J}B_{r'}(x_{k'})}f_{1,u}dV_g\phi(x_k)\rr\rm\\
\nn&\le&\|\n\phi\|_{L^\ity(\Si)}\sum_{x_k\in\mc J}\int_{B_{r'}(x_k)\bs\Cup_{k'=1}^kB_{r'}(x_{k'})}f_{1,u}d(\cd,x_k)dV_g\\
\nn&+&\|\phi\|_{L^\ity(\Si)}\int_{\Si\bs\Cup_{x_{k'}\in\mc J}B_{r'}(x_{k'})}f_{1,u}dV_g\\
\nn&<&r'\|\n\phi\|_{L^\ity(\Si)}\int_{\Cup_{x_{k'}\in\mc J}B_{r'}(x_{k'})}f_{1,u}dV_g\\
\nn&+&\e'\|\phi\|_{L^\ity(\Si)}\\
\nn&\le&r'\|\n\phi\|_{L^\ity(\Si)}+\e'\|\phi\|_{L^\ity(\Si)}\\
\lab{depsin}&\le&\fr{3}4\e\|\phi\|_{Lip(\Si)}.
\eea
Therefore, from $\dis{\eqr{depsout}}$ and $\dis{\eqr{depsin}}$ we deduce $\dis{\eqr{depsphi}}$, hence $\dis{\eqr{depslip}}$.
\epf\

\sec{Mapping sub-levels into weighted barycenters}\

The results we obtained about the sub-levels of $\dis{J_\rho}$ will be used in this section to build a map on $\dis{J_\rho^{-L}}$ which can be combined with the map $\dis{\Phi_\l}$ defined in Theorem $\dis{\ref{test}}$ to get a homotopy equivalence.\\
Precisely, we will get the following result:\\

\bth\lab{phipsi}$\dis{}$\\
Suppose $\dis{\rho\nin\L}$. Then, for $\dis{L}$ large enough there exist two continuous maps
$$\Phi:\g_{\st,\rho,\wt\a}\to J_\rho^{-L},\q\q\q\Psi:J_\rho^{-L}\to\g_{\st,\rho,\wt\a}$$
such that the $\dis{\Psi\c\Phi}$ is homotopically equivalent to $\dis{\mr{Id}_{\g_{\st,\rho,\wt\a}}}$.
\eth\

For the map $\dis{\Phi}$ we will choose $\dis{\Phi_\l}$ with a suitable $\dis{\l\gg0}$, whereas $\dis{\Psi}$ will be modeled on the retraction $\dis{\psi_{\rho,\a}}$ defined in Lemma $\dis{\ref{ret}}$.\\
The main issue is choosing properly the parameter $\dis{t\in[0,1]}$ in the join, so that it equals $\dis{0}$ or $\dis{1}$ as long as only one retraction is defined. The next lemma will give an estimate on the distance between the components of $\dis{\ph_{\l,\z}}$ and the respective weighted barycenter sets, thus giving a hint on when each retraction can or cannot be applied.\\

\ble\lab{dbar}$\dis{}$\\
Let $\dis{\s_i,\z,\ph_{\l,\z},\b_{ik}}$ be as in Theorem $\dis{\ref{test}}$. Then, for any $\dis{i=1,2}$, $\dis{\l>0}$, $\dis{\z\in\g_{\st,\rho,\wt\a}}$ one has
$$\fr{1}CF(\s_i,\l_{i,t})\le d_{Lip'}\lr f_{i,\ph_{\l,\z}},\Si_{\rho_i,\wt\a_i}\rr\le CF(\s_i,\l_{i,t})$$
with
$$F(\s_i,\l_{i,t}):=\sum_{x_{ik}\in\mc J_i}\fr{t_{ik}}{\max\{1,\l_{i,t}\}^{\min\lb2,\fr{1}{1+\b_{ik}}\rb}}$$
\ele\

\bpf$\dis{}$\\
It clearly suffices to give the proof for $\dis{i=1}$ and for large $\dis{\l_{1,t}}$.\\
For the upper bound, we will show that
$$d_{Lip'}\lr f_{1,\ph_{\l,\z}},\s_{1,\l}\rr\le C\sum_{x_{1k}\in\mc J_1}\fr{t_{1k}}{\l_{1,t}^{\min\lb2,\fr{1}{1+\b_{1k}}\rb}}$$
with
$$\s_{1,\l}:=\sum_{x_{1k}\in\mc J_1}t_{1k,\l}\d_{x_{1k}},\q\q\q t_{1k,\l}=t_{1k}\fr{\int_\Si\fr{\wt h_1}{\lr1+\l_{1,t}^2d(\cd,x_{1k})^{2\lr1+\b_{1k}\rr}\rr^2}\lr\sum_{x_{2k'}\in\mc J_2}\fr{t_{2k'}}{\lr1+\l_{2,t}^2d(\cd,x_{2k'})^{2\lr1+\b_{2k'}\rr}\rr^2}\rr^{-\fr{1}2}dV_g.}{\int_\Si\wt h_1e^{\ph_{1,\l,\z}}dV_g}.$$
From Lemma $\dis{\ref{log}}$, given any $\dis{\phi\in Lip(\Si)}$ with $\dis{\|\phi\|_{Lip(\Si)}\le1}$ we find
\bey
\lm\int_\Si\lr f_{1,\ph_{\l,\z}}-\s_{1,\l}\rr\phi dV_g\rm&=&\fr{1}{\int_\Si\wt h_1e^{\ph_{1,\l,\z}}dV_g}\lr\int_\Si\lr\wt h_1e^{\ph_{1,\l,\z}}-\int_\Si\wt h_1e^{\ph_{1,\l,\z}}dV_g\s_{1,\l}\rr\phi dV_g\rr\\
&\le&\fr{\l_{1,t}^2}{\l_{2,t}^2}\lm\int_\Si\wt h_1e^{\ph_{1,\l,\z}}\lr\phi-\sum_{x_{1k}\in\mc J_1}t_{1k,\l}\phi(x_{1k})\rr dV_g\rm\\
&=&\fr{\l_{1,t}^2}{\l_{2,t}^2}\lm\int_\Si\lr\sum_{x_{1k}\in\mc J_1}t_{1k}\fr{\fr{\wt h_1}{\lr1+\l_{1,t}^2d(\cd,x_{1k})^{2(1+\b_{1k'})}\rr^2}}{\sqrt{\sum_{x_{2k'}\in\mc J_2}\fr{t_{2k'}}{\lr1+\l_{2,t}^2d(\cd,x_{2k'})^{2\lr1+\b_{2k'}\rr}\rr^2}}}\rr(\phi-\phi(x_{1k}))dV_g\rm\\
&\le&\l_{1,t}^2\lm\int_\Si\sum_kt_{1k}\fr{\wt h_1}{\lr1+\l_{1,t}^2d(\cd,x_{1k})^{2(1+\b_{1k'})}\rr^2}(\phi-\phi(x_{1k}))dV_g\rm\\
&\le&\l_{1,t}^2\sum_kt_{1k}\int_\Si\fr{\wt h_1d(\cd,x_{1k})}{\lr1+\l_{1,t}^2d(\cd,x_{1k})^{2(1+\b_{1k'})}\rr^2}dV_g,
\eey
hence the estimate will follow if we show
$$\l_{1,t}^2\int_\Si\fr{\wt h_1d(\cd,x)}{\lr1+\l_{1,t}^2d(\cd,x)^{2(1+\b)}\rr^2}dV_g\le\fr{C}{\l_{1,t}^{\min\lb2,\fr{1}{1+\b}\rb}}$$
for any $\dis{x=x_{1k}}$, $\dis{\b=\b_{1k}}$.\\
We easily find
$$\l_{1,t}^2\int_{\Si\bs B_\d(x)}\fr{\wt h_1d(\cd,x)}{\lr1+\l_{1,t}^2d(\cd,x)^{2(1+\b)}\rr^2}dV_g\le\fr{C}{\l_{1,t}^2};$$
on the other hand, using normal coordinates and a change of variable we find, if $\dis{\b=0}$,
$$\l_{1,t}^2\int_{B_\d(x)}\fr{\wt h_1d(\cd,x)}{\lr1+\l_{1,t}^2d(\cd,x)^{2(1+\b)}\rr^2}dV_g\le\fr{C}{\l_{1,t}}\int_{B_{\l_{1,t}\d}(0)}\fr{|y|}{\lr1+|y|^2\rr^2}dy\le\fr{C}{\l_{1,t}};$$
if $\dis{x}$ is close to a point $\dis{p}$ with a singularity $\dis{\a}$, then
$$\l_{1,t}^2\int_{B_\d(x)}\fr{\wt h_1d(\cd,x)}{\lr1+\l_{1,t}^2d(\cd,x)^{2(1+\b)}\rr^2}dV_g\le\fr{C}{\l_{1,t}^\fr{1}{1+\b}}\int_{B_{\l_{i,t}^\fr{1}{1+\b}\d}(0)}\fr{\lm\l_{1,t}^\fr{\b-\a}{(1+\b)\a}y-\l_{1,t}^\fr{\b}{(1+\b)\a}p\rm^{2\a}|y|}{\lr1+|y|^{2(1+\b)}\rr^2}dy\le\fr{C}{\l_{1,t}^\fr{1}{1+\b}},$$
since the last integral is uniformly bounded (see \ci{carmal}, Proposition $\dis{4.1}$).
To give a lower bound, it suffices to prove that, however we take $\dis{\s=\s_\l}$, there exists a $\dis{1-Lip}$ function $\dis{\phi_\s}$ which satisfies
$$\lm\int_\Si\lr f_{1,\ph_{\l,\z}}-\s\rr\phi_\s dV_g\rm\ge\fr{1}C\sum_{x_{1k}\in\mc J_1}\fr{t_{1k}}{\l_{1,t}^{\min\lb2,\fr{1}{1+\b_{1k}}\rb}}.$$
Precisely, we choose
$$\phi_\s=\min_{x_{k'}\in\mc J'}d(\cd,x_{k'})\q\q\q\tx{if }\s=\sum_{x_{k'}\in\mc J'}t_{k'}\d_{x_{k'}}.$$
It holds
\bey
\lm\int_\Si\lr f_{1,\ph_{\l,\z}}-\s\rr\phi_\s dV_g\rm&=&\fr{1}{\int_\Si\wt h_1e^{\ph_{1,\l,\z}}dV_g}\lm\int_\Si\lr\wt h_1e^{\ph_{1,\l,\z}}-\int_\Si\wt h_1e^{\ph_{1,\l,\z}}dV_g\s\rr\phi_\s dV_g\rm\\
&=&\fr{1}{\int_\Si\wt h_1e^{\ph_{1,\l,\z}}dV_g}\int_\Si\wt h_1e^{\ph_{1,\l,\z}}\min_{k'}d(\cd,x_{k'})dV_g\\
&\ge&\fr{\l_{1,t}^2}{\l_{2,t}^2}\int_{\Si}\wt h_1e^{\ph_{1,\l,\z}}\min_{k'}d(\cd,x_{k'})dV_g\\
&\ge&\l_{1,t}^2\sum_{x_k\in\mc J_1}t_{1k}\int_{\Si}\fr{\wt h_1\min_{k'}d(\cd,x_{k'})}{\lr1+d(\cd,x_{1k})^{2(1+\b_{1k'})}\rr^2}dV_g.
\eey
Again, it is easy to see that any single integral outside a ball $\dis{B_\d(x_{1,k})}$ is greater or equal to constant times $\dis{\l_{1,t}^{-2}}$, since the number of $\dis{k'}$ is at most $\dis{K=K(\rho_1,\wt\a_1)}$. Therefore, it will suffice to show that any integral on the same ball can be estimated from below with constant times $\dis{\l_{1,t}^{-\fr{1}{1+\b_{1k}}}}$. Arguing as before,
$$\l_{1,t}^2\int_{B_\d(x)}\fr{\wt h_1\min_{k'}d(\cd,x_{k'})}{\lr1+\l_{1,t}^2d(\cd,x)^{2(1+\b)}\rr^2}dV_g\ge\fr{1}{C\l_{1,t}^\fr{1}{1+\b}}\int_{B_{\l_{i,t}^\fr{1}{1+\b}\d}(0)}\fr{\lm\l_{1,t}^\fr{\b-\a}{(1+\b)\a}y-\l_{1,t}^\fr{\b}{(1+\b)\a}p\rm^{2\a}\min_{k'}\lm y-\l_{1,t}^\fr{1}{1+\b}x_{k'}\rm}{\lr1+|y|^{2(1+\b)}\rr^2}dy.$$
To see that the last integral is bounded from above, we restrict ourselves to a portion of a ball where the minimum is attained by $\dis{x'=x_{k'}}$. Since the number of $\dis{k'}$ is uniformly bounded, for at least one index the portion we are considering measures at least $\dis{\fr{1}K}$ of the measure of the whole ball.\\
If we take $\dis{x'=x'_\l}$ so that $\dis{\l_{1t}^\fr{1}{1+\b}x'_\l}$ goes to infinity, the integral will tend to $\dis{+\ity}$ as well; if instead the last quantity converges, we will get the integral of a function which is uniformly bounded from both above and below, as in the proof of the upper estimates a few lines before.\\
The same argument works when we have no singularities in $\dis{B_\d(x)}$.
\epf\

From this lemma we also deduce a useful corollary:\\

\bco$\dis{}$\\
Let $\dis{\s_i,\z,\ph_{\l,\z}}$ be as in Theorem $\dis{\ref{test}}$ and $\dis{t_{ik,\l}}$ be as in Lemma $\dis{\ref{dbar}}$. Then, if $\dis{t\ne1}$ we have
$$f_{1,\ph_{\l,\z}}\us{\l\to\ity}\wk\wt\s_1:=\sum_{x_{1k}\in\mc J}\wt t_{1k}\d_{x_{1k}}$$
and if $\dis{t\ne0}$ we have
$$f_{2,\ph_{\l,\z}}\us{\l\to\ity}\wk\wt\s_2:=\sum_{x_{2k}\in\mc J}\wt t_{2k}\d_{x_{2k}}$$
where $\dis{\wt t_{ik}}$ verifies
$$\wt t_{ik}=\lim_{\l\to\ity}t_{ik,\l},\q\q\q\fr{t_{ik}}C\le\wt t_{ik}\le Ct_{ik}$$
\eco\

We are now in a position to prove Theorem $\dis{\ref{phipsi}}$.\\

\bpf[Proof of Theorem $\dis{\ref{phipsi}}$]$\dis{}$\\
Fix $\dis{C}$ as in Lemma $\dis{\ref{dbar}}$, $\dis{\e_0}$ as in Lemma $\dis{\ref{ret}}$ and apply Theorem $\dis{\ref{deps}}$ with $\dis{\e:=\fr{\e_0}{C^2}}$.\\
Take now $\dis{L=L(\e)>0}$ as in Theorem $\dis{\ref{deps}}$ and define $\dis{\Phi:=\Phi_{\l_0}}$ with $\dis{\l_0}$ such that $\dis{J(\ph_{\l,\z})\le-L}$ for any $\dis{\z\in\g_{\st,\rho,\wt\a}}$.\\
Define now, for $\dis{u\in J_\rho^{-L}}$,
\beq\lab{tt}
\wt t=\wt t(d_1,d_2):=\lb\bll0&\tx{if }d_2\ge\e\\\fr{\e-d_2}{2\e-d_1-d_2}&\tx{if }d_1,d_2<\e\\1&\tx{if }d_1\ge\e\ea\ry\q\q\q\tx{where }d_i=d_{Lip'}(f_{i,u},\Si_{\rho_i,\wt\a_i}).
\eeq
This quantity is always well-defined and continuous, because on $\dis{J_\rho^{-L}}$ at least one of $\dis{d_1}$ and $\dis{d_2}$ is less than $\dis{\e}$.\\
Consider now $\dis{\psi_i:=\psi_{\rho_i,\wt\a_i}}$ as in Lemma $\dis{\ref{ret}}$ and the push-forward $\dis{(\Pi_i)_*}$ of the maps $\dis{\Pi_i=\Si\to\g_i}$. We can now define the map $\dis{\Psi}$, from $\dis{J_\rho^{-L}}$ to $\dis{\g_{\st,\rho,\wt\a}}$:
$$\Psi(u)=\lr1-\wt t\rr(\Pi_1)_*(\psi_1(f_{1,u}))+\wt t(\Pi_2)_*(\psi_2(f_{2,u})).$$
This map is well-defined as well because, from the construction of $\dis{\wt t}$, when $\dis{\psi_1}$ is not defined one has $\dis{d_1\ge\e_0>\e}$, hence $\dis{\wt t=1}$, and similarly $\dis{\wt t=0}$ when $\dis{\psi_2}$ is not defined.\\
Let us now compose the maps $\dis{\Phi}$ and $\dis{\Psi}$ and see what happens if we let $\dis{\l}$ tend to $\dis{+\ity}$. From the previous corollary, $\dis{f_{i,\ph_{\l,\z}}}$ converges weakly to a barycenter $\dis{\wt\s_i}$ centered at the same points as $\dis{\s_i}$, and the same convergence still holds after applying $\dis{\psi_i}$ and $\dis{(\Pi_i)_*}$, since both are retractions. However, the coefficients in $\dis{\wt\s_i}$ are different from the ones in $\dis{\s_i}$, and moreover the parameter $\dis{t}$ in the join will be different in general from $\dis{\wt t}$.\\
Following this considerations, we will construct the homotopy between $\dis{\Psi\c\Phi}$ and the identity in three steps: first letting $\dis{\l}$ to $\dis{+\ity}$, then rescaling the coefficients in $\dis{\wt\s_i}$ and finally rescaling the parameter $\dis{t}$ in the join. Precisely, the homotopy map $\dis{H:\g_{\st,\rho,\wt\a}\x[0,1]\to\g_{\st,\rho,\wt\a}}$ will be the composition $\dis{H:=H_3*(H_2*H_1)}$, where:
$$H_1:(\z,s)=((1-t)\s_1+t\s_2,s)\to\lr1-\wt t\rr(\Pi_1)_*\lr\psi_1\lr f_{1,\ph_{\fr{\l_0}{1-s},\z}}\rr\rr+\wt t(\Pi_2)_*\lr\psi_2\lr f_{2,\ph_{\fr{\l_0}{1-s},\z}}\rr\rr$$
$$H_2:\lr\lr1-\wt t\rr\wt\s_1+\wt t\wt\s_2,s\rr\to\lr\lr1-\wt t\rr((1-s)\wt\s_1+s\s_1)+\wt t((1-s)\wt\s_2+s\s_2)\rr$$
$$H_3:\lr\lr1-\wt t\rr\s_1+\wt t\s_2,s\rr\to\lr\lr1-\lr(1-s)\wt t+st\rr\rr\s_1+\lr(1-s)\wt t+st\rr\s_2\rr.$$
Let us now verify that the maps are well defined.\\
In the definition of the map $\dis{H_1}$, Lemma $\dis{\ref{dbar}}$ ensures that the retraction $\dis{\psi_1}$ is defined if we have $\dis{F\lr\s_i,\fr{\l_{0i,t}}{1-s}\rr<\fr{\e_0}C}$. If the latter quantity is greater or equal to $\dis{\fr{\e_0}C}$, then $\dis{\psi_1}$ might not be defined, but in this case we have
$$d_1\ge d_{Lip'}\lr f_{1,\ph_{\fr{\l_0}{1-s},\z}},\Si_{\rho_1,\wt\a_1}\rr\ge\fr{\e_0}{C^2}=\e,$$
hence $\dis{\wt t=1}$ and therefore everything makes sense. For the same reason we can compose $\dis{\psi_2}$ and $\dis{\wt t}$.\\
In $\dis{H_2}$, the convex combination of $\dis{\wt\s_i}$ and $\dis{\s_i}$ are allowed in $\dis{(\g_i)_{\rho_i,\wt\a_i}}$ because the centers of the Dirac masses which define them are the same.\\
Finally, it is immediate to see that the composition makes sense, namely $\dis{H_1(\cd,1)=H_2(\cd,0)}$ and $\dis{H_2(\cd,1)=H_3(\cd,0)}$, that $\dis{H(\cd,0)=\Psi\circ\Phi}$ and $\dis{H(\cd,1)=\mr{Id}_{\g_{\st,\rho,\wt\a}}}$, and that everything is continuous.
\epf\

The existence of this homotopy map gives, through the functorial properties of homology, a simple but very important corollary:\\

\bco\lab{hom} Suppose $\dis{\rho\nin\L}$. Then, for $\dis{L}$ large enough the map $\dis{\Phi}$ defined in theorem $\dis{\ref{phipsi}}$ induces an immersion of homology groups
$$H_q(\g_{\st,\rho,\wt\a})\sr{\Phi_{*,q}}\inc H_q\lr J_\rho^{-L}\rr\q\q\q\fa\,q\in\N.$$
Therefore, in particular, $\dis{\wt\b_q\lr J_\rho^{-L}\rr\ge\wt\b_q(\g_{\st,\rho,\wt\a})}$ for any $\dis{q\in\N}$.
\eco\

\sec{The weighted barycenter sets}\

In this section, we will provide information about the topology and the homology of the space $\dis{\g_{\st,\rho,\wt\a}=(\g_1)_{\rho_1,\wt\a_1}\st(\g_2)_{\rho_2,\wt\a_2}}$, whose importance in the study of the problem arose clearly in Theorem $\dis{\ref{phipsi}}$ and Corollary $\dis{\ref{hom}}$.\\
First of all, we notice that most information can be deduced by studying the weighted barycenters spaces $\dis{(\g_i)_{\rho_i,\wt\a_i}}$. Proposition $\dis{\ref{join}}$ shows how the homology groups of the join depend on the ones of the spaces which form it.\\
Some of the results contained in this section will be inspired by \ci{car}, where weighted barycenters centered at $\dis{\Si}$ are studied.\\
Moreover, it is easy to see that if one of the two spaces is contractible, then the join is contractible as well. In fact, if $\dis{H}$ is a homotopy equivalence between $\dis{X}$ and a point, then
$$((1-t)x+ty,s)\to(1-t)H(x,s)+ty$$
is a homotopy equivalence between $\dis{X\st Y}$ and the cone based in $\dis{Y}$, which is contractible.\\
Therefore, it suffices to restrict our study to the weighted barycenter sets $\dis{(\g_i)_{\rho_i,\wt\a_i}}$. In the following, we will omit the indices $\dis{i=1,2}$ and consider a generic weighted barycenters set $\dis{\g_{\rho,\a}}$ with the multi-indices $\dis{\a=(\a_1,\ds,\a_L)}$ such that $\dis{\a_l\le\a_{l+1}}$ and singular points $\dis{p_1,\ds,p_L}$ satisfy $\dis{\chi_\a(p_l)=1+\a_l<1}$.\\
To start with, following \ci{car} we consider $\dis{\g_{\rho,\a}}$ as a union of strata of the kind
$$\g^{K,\mc I}=\lb\sum_{k=1}^Kt_k\d_{q_k}+\sum_{l\in\mc I}s_l\d_{p_l};\,q_k\in\Si,\,t_k\ge0,\,s_l\ge0,\,\sum_{k=1}^Kt_k+\sum_{l\in\mc I}s_l=1\rb\q\tx{for }K\in\N\cup\{0\},\,\mc I\sub\{1,\ds,L\}.$$
One can easily notice that each of these strata is a union of manifolds whose maximal dimension is $\dis{2K+|\mc I|-1}$. Considering only the strata which are maximal with respect to the inclusion, we write a unique decomposition
\beq\lab{maxstr}
\g_{\rho,\a}=\Cup_{h=1}^H\g^{K_h,\mc I_h}.
\eeq
It is easy to see how the strata depend on the position of $\dis{\rho}$ with respect to the $\dis{\a_l}$'s. A stratum $\dis{\g^{K,\mc I}}$ is contained in $\dis{\g_{\rho,\a}}$ if and only if
\beq\lab{str}
\rho>4\pi\lr K+\sum_{l\in\mc I}(1+\a_l)\rr.
\eeq
Moreover, we notice that a stratum $\dis{\g^{K,\mc I}}$ is contained in $\dis{\g^{K',\mc I'}}$ if and only if $\dis{|\mc I\bs\mc I'|\le K'-K}$. Therefore, the maximality of an existing stratum is equivalent to the condition
$$\rho\le4\pi\min\lb K+1+\sum_{l\in\mc I\bs\{\max\mc I\}}(1+\a_l),K+\sum_{l\in\mc I\cup\{\min(\{1,\ds,L\}\bs\mc I)\}}(1+\a_l)\rb,$$
and the equality sign is excluded if we take $\dis{\rho\nin\L}$.\\
Notice that in the regular case the decomposition in maximal strata is just $\dis{\g_{\rho,\es}=\g^{K,\es}=\g^K}$, with $\dis{K}$ such that $\dis{\rho\in(4K\pi,4(K+1)\pi)}$, and all the strata are of the kind $\dis{\g^{K',\es}=\g^{K'}}$ for $\dis{K'=1,\ds,K}$. However, in the regular case Proposition $\dis{\ref{homb}}$ gives already full information about homology of the barycenters.\\
In the general case the decomposition in strata makes more difficult the computation of the homology groups. Nonetheless, we can still obtain information on the homology of $\dis{\g_{\rho,\a}}$ with an estimate from below of its Betti numbers.\\

\bth\lab{homwb}$\dis{}$\\
Suppose $\dis{\g_{\rho,\a}}$ has the following decomposition in maximal strata:
\beq\lab{maxstr1}
\g_{\rho,\a}=\Cup_{h=1}^H\g^{K_h,\mc I_h}\cup\Cup_{h'=1}^{H'}\g^{K'_{h'},\mc I'_{h'}},
\eeq
with $\dis{1\nin\mc I_h}$ for any $\dis{h=1,\ds,H}$. Then,
$$\wt\b_q(\g_{\rho,\a})\ge\sum_{h=1}^H\bin{K_h+|\mc I_h|+\ls\fr{-\chi(\Si)}2\rs}{|\mc I_h|+\ls\fr{-\chi(\Si)}2\rs}\d_{q,2K_h+\mc I_h-1}.$$
In particular, if $\dis{h\ge1}$, then $\dis{\wt\b_q(\g_{\rho,\a})\ne0}$ for some $\dis{q\ne0}$.
\eth\

We will start by seeing the cases which are not covered by the previous theorem, that is when every maximal stratum is defined by a multi-index containing the index $\dis{1}$.\\
In this case, we find out that $\dis{\g_{\rho,\a}}$ is contractible, so in conclusion we get a necessary and sufficient condition for the contractibility of $\dis{\g_{\rho,\a}}$.

\ble$\dis{}$\\
Suppose $\dis{\g_{\rho,\a}}$ has the decomposition $\dis{\eqr{maxstr}}$ in maximal strata, with $\dis{p_1,\ds,p_L}$ such that $\dis{\a_1\le\ds\le\a_L}$. Then, the following conditions are equivalent:
\ben
\it $\dis{\g_{\rho,\a}}$ is star-shaped with respect to $\dis{\d_{p_1}}$.
\it There exists some $\dis{l\in\{1,\ds,L\}}$ such that $\dis{\g_{\rho,\a}}$ is star-shaped with respect to $\dis{\d_{p_l}}$.
\it $\dis{\g^{K_h,\mc I_h}}$ is star-shaped with respect to $\dis{\d_{p_1}}$ for any $\dis{h\in H}$.
\it There exists some $\dis{l\in\{1,\ds,L\}}$ such that $\dis{\g^{K_h,\mc I_h}}$ is star-shaped with respect to $\dis{\d_{p_l}}$ for any $\dis{h\in\{1,\ds,H\}}$.
\een
Moreover, each of these conditions implies that $\dis{\g_{\rho,\a}}$ is contractible.
\ele\ 

\bpf$\dis{}$\\
The contractibility of $\dis{\g_{\rho,\a}}$ follows trivially from its star-shapedness, so it suffices to prove the equivalences between the conditions.\\
The following implications are evident:
$$\mi1\To\mi2,\q\q\q\mi3\To\mi1,\q\q\q\mi3\To\mi4,\q\q\q\mi4\To\mi2;$$
therefore, we suffice to show that $\dis{\mi2}$ implies $\dis{\mi1}$ and $\dis{\mi1}$ implies $\dis{\mi3}$.\\
We will start by showing that if $\dis{\g_{\rho,\a}}$ is star-shaped with respect to some $\dis{p_{\wt l}}$, then the same holds with $\dis{p_1}$.\\
We notice immediately that star-shapedness of $\dis{\g_{\rho,\a}}$ is equivalent to saying that for any stratum $\dis{\g^{K,\mc I}\sub\g_{\rho,\a}}$ we have $\dis{\g^{K,\mc I\cup\wt l}\sub\g_{\rho,\a}}$; moreover, we recall that the existence of a stratum within $\dis{\g_{\rho,\a}}$ means $\dis{\eqr{str}}$. Let us now suppose condition $\dis{\mi 2}$ occurs for $\dis{\wt l>1}$, that is
$$\rho>4\pi\lr K+\sum_{l\in\mc I}(1+\a_l)\rr\q\q\q\To\q\q\q\rho>4\pi\lr K+\sum_{l\in\mc I\cup\lb\wt l\rb}(1+\a_l)\rr,$$
and let us recall that we are assuming $\dis{\a_l\le\a_{l+1}}$ for any $\dis{l}$. This implies
$$\rho>4\pi\lr K+\sum_{l\in\mc I\cup\lb\wt l\rb}(1+\a_l)\rr\ge4\pi\lr K+\sum_{l\in\mc I\cup\{1\}}(1+\a_l)\rr,$$
that is star-shapedness of $\dis{\g_{\rho,\a}}$ with respect to $\dis{p_1}$.\\
Suppose now, by contradiction, that condition $\dis{\mi3}$ holds but condition $\dis{\mi1}$ does not, that is $\dis{\g_{\rho,\a}}$ is star-shaped with respect to $\dis{p_1}$ but it contains a maximal stratum $\dis{\g^{K,\mc I}}$ which is not.\\
Then, star-shapedness of $\dis{\g_{\rho,\a}}$ with respect to $\dis{\d_{p_1}}$ implies the existence of a stratum $\dis{\g^{K,\mc I\cup\{1\}}\sub\g_{\rho,\a}}$, which contains properly $\dis{\g^{K,\mc I}}$, thus contradicting its maximality.
\epf\

Let us now see what happens if we are in a scenario which is opposite to the previous lemma, that is some index $\dis{j}$ is not contained in any multi-index which defines the strata.\\
The following lemma shows that this situation produces some non-trivial homology.\\

\ble\lab{nontriv}$\dis{}$\\
Suppose $\dis{K\in\N}$, $\dis{\mc I\sub\{1,\ds,L\}}$ and $\dis{\wt l\nin\mc I}$ and define
$$\wt\g^{K,\mc I,\wt l}:=\Cup_{\mc I'\sub\mc I\cup\lb\wt l\rb,\,|\mc I'|=|\mc I|}\g^{K,\mc I}.$$
Then, it holds
$$\wt H_q\lr\wt\g^{K,\mc I,\wt l}\rr=\lb\bll\Z^{\bin{K+|\mc I|+\ls\fr{-\chi(\Si)}2\rs}{|\mc I|+\ls\fr{-\chi(\Si)}2\rs}}&\tx{if }q=2K+|\mc I|-1\\0&\tx{if }q\ne2K+|\mc I|-1\ea\ry.$$
\ele\

The proof of the lemma will use the Mayer-Vietoris exact sequence.\\
Actually, when applying the Mayer-Vietoris sequence the sets $\dis{A}$ and $\dis{B}$ should be open. If they are not, we are implicitly considering two suitable open neighborhoods in their stead.\\
The existence of such neighborhoods follows from the properties of the weighted barycenters, which can be deduced by arguing as in \cite{car}, Section $\dis{2}$ and \cite{carmal}, Section $\dis{3}$.

\bpf$\dis{}$\\
We proceed by double induction on $\dis{K}$ and $\dis{|\mc I|}$.\\
If $\dis{\mc I=\es}$ we have $\dis{\wt\g^{K,\es,l}=\g^{K,\es}=\g^K}$ so the claim follows by Proposition $\dis{\ref{homb}}$.\\
If $\dis{K=0}$, any stratum $\dis{\g^{0,\mc I'}}$ is actually the $\dis{(|\mc I'|-1)}$-simplex $\dis{\ls\d_{p_{l_1}},\ds,\d_{p_{l_{|\mc I'|}}}\rs}$ if we can write $\dis{\mc I'=\lb l_1,\ds,l_{|\mc I'|}\rb}$. Therefore, $\dis{\wt\g^{0,\mc I,\wt l}}$ is the boundary of the $\dis{|\mc I|}$-simplex with vertices in $\dis{\d_{p_l}}$ for $\dis{l\in\mc I\cup\lb\wt l\rb}$; hence, it is homeomorphic to the sphere $\dis{\S^{|\mc I|-1}}$ and the claim follows also in this case.\\
Suppose now that the lemma is true for $\dis{K-1,\mc I}$ and for any $\dis{K,\mc I_0}$ with $\dis{\lm\mc I_0\rm=|\mc I|-1}$.\\
Being $\dis{\g^{K,\mc I,\wt l}}$ union of manifolds of dimension less or equal to $\dis{2K+|\mc I|-1}$, all the higher homology groups are trivial.\\
To compute the other groups, we write $\dis{\wt\g^{K,\mc I,\wt l}=A\cup B}$ with
$$A=\g^{K,\mc I},\q\q\q B=\Cup_{l\in\mc I}\g^{K,\mc I\bs l\cup\lb\wt l\rb}$$
and consider the Mayer-Vietoris sequence. The set $\dis{B}$ is star-shaped with respect to $\dis{\d_{p_{\wt l}}}$ whereas $\dis{A}$ is star-shaped with respect to $\dis{\d_{p_l}}$ for any $\dis{l\in\mc I}$, hence we can write
$$0=\wt H_q(A)\pl\wt H_q(B)\to\wt H_q(A\cup B)\to\wt H_{q-1}(A\cap B)\to\wt H_{q-1}(A)\pl\wt H_{q-1}(B)=0,$$
that is $\dis{\wt H_q(A\cup B)=\wt H_{q-1}(A\cap B)}$. Moreover, this set can be written as
$$A\cap B=C\cup D,\q\q\q C:=\g^{K-1,\mc I\cup\lb\wt l\rb},\q\q\q D:=\Cup_{l\in\mc I}\g^{K,\mc I\bs\{l\}}.$$
As before, $\dis{C}$ is contractible, whereas we can write $\dis{D=\wt\g^{K,\mc I\bs\{l\},\{l\}}}$ for any $\dis{l\in\mc I}$ and $\dis{C\cap D=\wt\g^{K-1,\mc I,\wt l}}$. Therefore, by inductive hypothesis we know the homology of these sets and we can apply again Mayer-Vietoris. If $\dis{q<2K+|\mc I|-1}$ we get
$$0=\wt H_{q-1}(C)\pl\wt H_{q-1}(D)\to\wt H_{q-1}(C\cup D)\to\wt H_{q-2}(C\cap D)\to\wt H_{q-2}(C)\pl\wt H_{q-2}(D)=0,$$
that is
$$\wt H_q(A\cup B)=\wt H_{q-1}(A\cap B)=\wt H_{q-1}(C\cup D)=\wt H_{q-2}(C\cap D)=0.$$
Finally, for the last homology group we get
$$0=\wt H_{2K+|\mc I|-2}(C\cap D)\to\wt H_{2K+|\mc I|-2}(C)\pl\wt H_{2K+|\mc I|-2}(D)\to\wt H_{2K+|\mc I|-2}(C\cup D)\to$$
$$\to\wt H_{2K+|\mc I|-3}(C\cap D)\to\wt H_{2K+|\mc I|-3}(C)\pl\wt H_{2K+|\mc I|-3}(D)=0.$$
Hence, by the inductive hypothesis and the properties of binomial coefficients,
\bey
\wt H_{2K+|\mc I|-1}(A\cup B)&=&\wt H_{2K+|\mc I|-2}(C\cup D)\\
&=&\wt H_{2K+|\mc I|-2}(D)\pl\wt H_{2K+|\mc I|-3}(C\cap D)\\
&=&\Z^{\bin{K+|\mc I|+\ls\fr{-\chi(\Si)}2\rs-1}{|\mc I|+\ls\fr{-\chi(\Si)}2\rs-1}}\pl\Z^{\bin{K+|\mc I|+\ls\fr{-\chi(\Si)}2\rs-1}{|\mc I|+\ls\fr{-\chi(\Si)}2\rs}}\\
&=&\Z^{\bin{K+|\mc I|+\ls\fr{-\chi(\Si)}2\rs-1}{|\mc I|+\ls\fr{-\chi(\Si)}2\rs-1}+\bin{K+|\mc I|+\ls\fr{-\chi(\Si)}2\rs-1}{|\mc I|+\ls\fr{-\chi(\Si)}2\rs}}\\
&=&\Z^{\bin{K+|\mc I|+\ls\fr{-\chi(\Si)}2\rs}{|\mc I|+\ls\fr{-\chi(\Si)}2\rs}},
\eey
which is what we wanted.
\epf\

Finally, we see how the sets defined in the previous lemma affect the homology of $\dis{\g_{\rho,\a}}$.\\

\bpf[Proof of Theorem $\dis{\ref{homwb}}$]$\dis{}$\\
We proceed by induction on $\dis{H}$. If $\dis{H=0}$ there is nothing to prove.\\
Suppose now the theorem holds true for $\dis{H-1}$. Then it also holds for $\dis{H}$ when $\dis{q\ne2K_H+|\mc I_H|-1}$.\\
For $\dis{q=2K_H+|\mc I_H|-1}$, we notice that $\dis{\wt\g^{K_H,\mc I_H,1}\sub\g_{\rho,\a}}$, since the coefficients $\dis{\a_l}$ are non-increasing; hence we can apply Mayer-Vietoris sequence by writing $\dis{\g_{\rho,\a}=A\cup B}$ with
$$A=\wt\g^{K_H,\mc I_H,1},\q\q\q B=\Cup_{h=1}^{H-1}\g^{K_h,\mc I_h}\cup\Cup_{h'=1}^{H'}\g^{K'_{h'},\mc I'_{h'}}.$$
By a dimensional argument we have $\dis{\wt H_{2K_H+|\mc I_H|-1}(A\cap B)=0}$, so we get
$$0=\wt H_{2K_H+|\mc I_H|-1}(A\cap B)\to\wt H_{2K_H+|\mc I_H|-1}(A)\pl\wt H_{2K_H+|\mc I_H|-1}(B)\to\wt H_{2K_H+|\mc I_H|-1}(A\cup B)\to\ds$$
which means, by the exactness of the Mayer-Vietoris sequence,
$$\wt H_{2K_H+|\mc I_H|-1}(A)\pl\wt H_{2K_H+|\mc I_H|-1}(B)\inc\wt H_{2K_H+|\mc I_H|-1}(A\cup B).$$
Therefore, applying the inductive hypothesis and Lemma $\dis{\ref{nontriv}}$,
we get
\bey
\wt\b_{2K_H+|\mc I_H|-1}(A\cup B)&\ge&\wt\b_{2K_H+|\mc I_H|-1}(A)+\wt\b_{2K_H+|\mc I_H|-1}(B)\\
&\ge&\bin{K+|\mc I|+\ls\fr{-\chi(\Si)}2\rs}{|\mc I|+\ls\fr{-\chi(\Si)}2\rs}\\
&+&\sum_{h=1}^{H-1}\bin{K_h+|\mc I_h|+\ls\fr{-\chi(\Si)}2\rs}{|\mc I_h|+\ls\fr{-\chi(\Si)}2\rs}\d_{2K_H+|\mc I_H|-1,2K_h+\mc I_h-1}\\
&=&\sum_{h=1}^H\bin{K_h+|\mc I_h|+\ls\fr{-\chi(\Si)}2\rs}{|\mc I_h|+\ls\fr{-\chi(\Si)}2\rs}\d_{2K_H+|\mc I_H|-1,2K_h+\mc I_h-1},
\eey
hence the claim.
\epf\

Finally, by Proposition $\dis{\ref{join}}$, we get information on the homology of the join.\\

\bco\lab{homjoin}$\dis{}$\\
Suppose $\dis{(\g_i)_{\rho_i,\wt\a_i}}$ has the decomposition $\dis{\eqr{maxstr1}}$ in maximal strata, with $\dis{H_i,K_1,\ds,K_{H_i}\in\N}$ and $\dis{\mc I_{i1},\ds,\mc I_{iH_i}\sub\{1,\ds,L_i\}}$. Then, it holds
$$\sum_{q=0}^{+\ity}\wt\b_q(\g_{\st,\rho,\wt\a})\ge\sum_{h_1=1}^{H_1}\sum_{h_2=1}^{H_2}\bin{K_{h_1}+|\mc I_{h_1}|+\ls\fr{-\chi(\Si)}2\rs}{|\mc I_{h_1}|+\ls\fr{-\chi(\Si)}2\rs}\bin{K_{h_2}+|\mc I_{h_2}|+\ls\fr{-\chi(\Si)}2\rs}{|\mc I_{h_2}|+\ls\fr{-\chi(\Si)}2\rs}.$$
In particular, if $\dis{h_1,h_2\ge1}$, then $\dis{\wt\b_q(\g_{\st,\rho,\wt\a})\ne0}$ for some $\dis{q\ne0}$.
\eco\

\sec{Examples and conclusion}\

Before proving the main results of this paper, let us see some examples of how they can be applied in dependence of the points $\dis{p_j}$ and the coefficients $\dis{\a_{ij}}$.
Since the condition $\dis{\eqr{rhoi}}$ only depends on the $\dis{\wt\a_{il}}$, for simplicity we will take all the $\dis{\a_{ij}}$ to be negative.\\

\bex$\dis{}$\\
Consider, for instance, the case of four singular points $\dis{p_{11}}$, $\dis{p_{12}}$, $\dis{p_{21}}$, $\dis{p_{22}}$, two of which having negative singularities $\dis{\wt\a_{11}=-\fr{3}4}$, $\dis{\wt\a_{12}=-\fr{1}3}$ for the first component, and the others having negative coefficients $\dis{\wt\a_{21}=-\fr{1}2}$, $\dis{\wt\a_{22}=-\fr{1}4}$ for the second component.\\
What are the possible values $\dis{\rho_1,\rho_2}$ which yield condition $\dis{\eqr{rhoi}}$? For ${\rho_1}$ we have to verify the condition with $\dis{\mc I_1=\es}$ and $\dis{\mc I_1=\{2\}}$, since we only have two points; therefore, $\dis{\rho_1}$ has to be either between $\dis{K}$ and $\dis{K+1+\wt\a_{11}=K+\fr{1}4}$ or between $\dis{K+1+\wt\a_{12}=K+\fr{2}3}$ and $\dis{K+2+\wt\a_{21}+\wt\a_{22}=K+\fr{11}{12}}$ for some integer $\dis{K}$. Hence, the admissible range for $\dis{\rho_1}$ is the union of the intervals
$$4\pi\lr0,\fr{1}4\rr\cup4\pi\lr\fr{2}3,\fr{11}{12}\rr\cup4\pi\lr1,1+\fr{1}4\rr\cup4\pi\lr1+\fr{2}3,1+\fr{11}{12}\rr\cup\ds.$$
Concerning $\dis{\rho_2}$, we can choose again $\dis{\mc I_2=\es}$ or $\dis{\mc I_2=\{2\}}$; as before, we get that $\dis{\rho_2}$ must lie in an interval of the kind $\dis{\lr K,K+\fr{1}2\rr}$ or $\dis{\lr K+\fr{3}4,K+\fr{5}4\rr}$ for some integer $\dis{K}$; anyway, the latter interval overlaps $\dis{\lr K+1,K+1+\fr{1}2\rr}$, so the range given by condition $\dis{\eqr{rhoi}}$ for $\dis{\rho_2}$ can be written as
$$4\pi\lr 0,\fr{1}2\rr\cup4\pi\lr\fr{3}4,\fr{3}2\rr\cup\pi\lr1+\fr{3}4,1+\fr{3}2\rr\cup\pi\lr2+\fr{3}4,2+\fr{3}2\rr\cup\ds.$$
\eex\

We are now in condition to finally prove the theorems stated in Section $\dis{1}$.\\

\bpf[Proof of Theorem $\dis{\ref{ex}}$]$\dis{}$\\
Suppose condition $\dis{\eqr{rhoi}}$ holds. This means that $\dis{(\g_i)_{\rho_i,\wt\a_i}}$ contains the stratum $\dis{\g^{K_i,\mc I_i}}$ and does not contain $\dis{(\g_i)^{K_i,\mc I_i\cup\{1\}}}$, for both $\dis{i=1,2}$. This stratum has to be contained in a maximal one $\dis{(\g_i)^{K'_i,\mc I'_i}}$ with $\dis{1\nin\mc I'_i}$, since otherwise we would have $\dis{(\g_i)^{K_i,\mc I_i\cup\{1\}}\sub(\g_i)^{K'_i,\mc I'_i}\sub(\g_i)_{\rho_i,\wt\a_i}}$.\\
Therefore, by Lemma $\dis{\ref{nontriv}}$, both $\dis{(\g_i)_{\rho_i,\wt\a_i}}$'s have non-trivial homology. Moreover, Corollary $\dis{\ref{homjoin}}$ ensures that $\dis{\g_{\st,\rho,\wt\a}}$ has non-trivial homology as well, and by Corollary $\dis{\ref{hom}}$ the same holds for $\dis{J_\rho^{-L}}$ if $\dis{L}$ is large enough.\\
Suppose now that the system $\dis{\eqr{toda}}$ has no solutions. Then, by Lemma $\dis{\ref{def}}$, $\dis{J_\rho^{-L}}$ should be a deformation retract of $\dis{J_\rho^L}$ for any $\dis{L}$. On the other hand, Corollary $\dis{\ref{contr}}$ says that for large $\dis{L}$ the sub-level $\dis{J_\rho^L}$ is contractible, whereas $\dis{J_\rho^{-L}}$ cannot be, having some non-trivial homology groups. Therefore, we are contradicting the assumption of having no solutions.
\epf\

\bpf[Proof of Theorem $\dis{\ref{mult}}$]$\dis{}$\\
Under the assumptions of the theorem, we can decompose each $\dis{(\g_i)_{\rho_i,\wt\a_i}}$ in maximal strata
$$(\g_i)_{\rho_i,\wt\a_i}=\Cup_{h_i=1}^{H_i}\g^{K_{h_i},\mc I_{h_i}}\cup\Cup_{h'_i=1}^{H'_i}\g^{K'_{h'_i},\mc I'_{h'_i}}.$$
Take $\dis{D}$ as the set of initial data such that $\dis{J_\rho}$ is a Morse function, which by Lemma $\dis{\ref{dense}}$ is a dense open set. Applying Lemma $\dis{\ref{num}}$ and Corollaries $\dis{\ref{hom}}$ and $\dis{\ref{homwb}}$ we get
\bey\#\tx{Solutions of }\eqr{toda}&\ge&\sum_{q=0}^{+\ity}\wt\b_q\lr J_\rho^{-L}\rr\\
&\ge&\sum_{q=0}^{+\ity}\wt\b_q(\g_{\st,\rho,\wt\a})\\
&\ge&\sum_{h_1,h_2}\bin{K_{h_1}+|\mc I_{h_1}|+\ls\fr{-\chi(\Si)}2\rs}{|\mc I_{h_1}|+\ls\fr{-\chi(\Si)}2\rs}\bin{K_{h_2}+|\mc I_{h_2}|+\ls\fr{-\chi(\Si)}2\rs}{|\mc I_{h_2}|+\ls\fr{-\chi(\Si)}2\rs},
\eey
that is the thesis of Theorem $\dis{\ref{mult}}$.
\epf\

\bpf[Proof of Theorem $\dis{\ref{exmult}}$]$\dis{}$\\
Due to the assumption $\dis{\rho_2<4\pi(1+\wh\a_2)}$, we can write
$$(\g_2)_{\rho_2,\wt\a_2}:=\lb\sum_{l\in\mc I}s_{2l}\d_{p_{2l}};\,s_l\ge0,\,\sum_{l\in\mc I}s_l=1,\,4\pi\sum_{l\in\mc I}(1+\wt\a_{2l})<\rho\rb.$$
If this set is not empty, that is if $\dis{\wh\a_2>\min_j\a_{2j}}$, we can still consider $\dis{\Phi_\l}$ as in Theorem $\dis{\ref{test}}$, since again, by construction, $\dis{d(\g_1,p_{2l})\ge\d>0}$ for any $\dis{l\in\{0,\ds,L_2\}}$.\\
Therefore we have, as in Theorem $\dis{\ref{test}}$, a map $\dis{\Phi:\g_{\st,\rho,\wt\a}\to J_\rho^{-L}}$ and, as in Theorem $\dis{\ref{phipsi}}$, $\dis{\Psi:J_\rho^{-L}\to\g_{\st,\rho,\wt\a}}$ such that $\dis{\Psi\c\Phi\seq\mr{Id}_{\g_{\st,\rho,\wt\a}}}$. Hence, the sublevels inherit the homology of the join, so existence and multiplicity of solutions follow by the estimating the Betti numbers as in Theorem $\dis{\ref{homwb}}$.\\

On the other hand, if $\dis{\wh\a_2=\min_j\a_{2j}}$, then the set $\dis{(\g_2)_{\rho_2,\wt\a_2}}$ is empty. However, $\dis{\Phi_\l}$ can still be defined on $\dis{(\g_1)_{\rho_1,\wt\a_1}}$ by restricting the map in Theorem $\dis{\ref{test}}$ to the end $\dis{t=0}$ of the join. Since we are just considering a restriction of the map, the estimates of the theorem still hold.\\
Moreover, being $\dis{\rho_2}$ small enough, Lemma $\dis{\ref{intbr}}$ can only hold for $\dis{i=1}$, so in Theorem $\dis{\ref{deps}}$ we must have $\dis{f_{1,u}}$ to be arbitrarily close to $\dis{\Si_{\rho_1,\wt\a_1}}$ as $\dis{J_\rho}$ is lower. Therefore, we can define $\dis{\Psi:J_\rho^{-L}\to(\g_1)_{\rho_1,\wt\a_1}}$ by $\dis{\Psi(u)=(\Pi_1)_*\psi_1(f_{1,u})}$ (with $\dis{\psi_1:=\psi_{\rho_1,\wt\a_1}}$ as in Lemma $\dis{\ref{ret}}$).\\
A homotopy map between $\dis{\Psi\c\Phi}$ and $\dis{\mr{Id}_{(\g_1)_{\rho_1,\wt\a_1}}}$ is given by restricting to $\dis{t=0}$ the map $\dis{H}$ defined in the proof of Theorem $\dis{\ref{phipsi}}$. Therefore, we can again deduce existence and multiplicity of solution by estimating the number of solutions as in Section $\dis{6}$.
\epf\

\bpf[Proof of Theorem $\dis{\ref{sph}}$]$\dis{}$\\
From the upper bound on $\dis{\rho_1}$ and $\dis{\rho_2}$, the elements in the barycenter sets $\dis{\Si_{\rho_i,\wt\a_i}}$ can be written in the simpler form
$$\Si_{\rho_i,\wt\a_i}:=\lb\sum_{l\in\mc I_i}s_{il}\d_{p_{il}};\,s_{il}\ge0,\,\sum_{l\in\mc I_i}s_{il}=1,\,4\pi\sum_{l\in\mc I_i}(1+\wt\a_{il})<\rho\rb.$$
We will assume both of these barycenter sets to be non-empty, since if one is empty we can modify the argument as in the proof of Theorem $\dis{\ref{exmult}}$ and if both are empty we are in the coercive case covered by Corollary $\dis{\ref{mt}}$.\\
As in Theorem $\dis{\ref{test}}$, we can build a map $\dis{\Phi_{\l}}$ from the join $\dis{\Si_{\st,\rho,\wt\a}=\Si_{\rho_1,\wt\a_1}\st\Si_{\rho_2,\wt\a_2}}$ by simply restricting the original $\dis{\Phi_\l}$ to these particular types of barycenters.\\
Since $\dis{d(p_{1l},p_{2l'})\ge\d>0}$ for any $\dis{l,l'}$, with the same argument we get $\dis{\Phi_\l\us{\l\to+\ity}\to-\ity}$ uniformly.\\
Moreover, from Theorem $\dis{\ref{deps}}$ (which also works when $\dis{g(\Si)=0}$), we get that $\dis{f_{i,u}}$ is $\dis{\e}$-close to $\dis{\Si_{\rho_i,\wt\a_i}}$ for some $\dis{i=1,2}$, so we can build a map $\dis{\Psi:J_\rho^{-L}\to\Si_{\st,\rho,\wt\a}}$ by taking $\dis{\wt t}$ as in $\dis{\eqr{tt}}$ and setting
$$\Psi(u):=\lr1-\wt t\rr\psi_1(f_{1,u})+\wt t\psi_2(f_{2,u}).$$
In the same way as in Theorem $\dis{\ref{phipsi}}$ we also prove the homotopy equivalence between $\dis{\Psi\circ\Phi}$ and $\dis{\mr{Id}_{\Si_{\st,\rho,\wt\a}}}$, so the sublevels inherit the homology of the join. We finally obtain the existence and multiplicity result by estimating its homology groups through Theorem $\dis{\ref{homwb}}$.\\
We can apply the latter theorem because, since regular points are not allowed in $\dis{\Si_{\rho_i,\wt\a_i}}$, this set coincides with $\dis{(\g_i)_{\rho_i,\wt\a_i}}$ for any simple closed curve $\dis{\g_i}$ (or also any subset of $\dis{\Si}$) which contains the points $\dis{\{p_{i1},\ds,p_{iL_i}\}}$.
\epf

\sec*{Acknowledgments}\

The author has been supported by the FIRB project \ti{Analysis and Beyond}, by the PRIN \ti{Variational and perturbative aspects of nonlinear differential problems} and by the Mathematics Department at the University of Warwick.\\
The author would like to thank Professor Andrea Malchiodi for the support and for the discussions concerning the topic of this paper.\\
Gratitude is also expressed to Professor Daniele Bartolucci for the discussions about the subject.

\bibliography{todasing4}
\bst{abbrv}

\edo